\documentclass[12pt]{article}

\usepackage{hyperref}
\hypersetup{
    unicode = false,
    pdftoolbar = true,
    pdfmenubar = true,
    pdffitwindow = true,
    pdftitle = {solar},
    pdfauthor = {Audet et al.},
    pdfsubject = {meta},
    pdfkeywords = {meta},
    pdfnewwindow = true,
    colorlinks = true,
    linkcolor = blue,
    citecolor = blue,
    filecolor = black,
    urlcolor = blue,
    breaklinks = true
}
\usepackage[margin=1in]{geometry}
\usepackage{booktabs}
\usepackage[utf8]{inputenc}
\usepackage{comment}
\usepackage[english]{babel}
\usepackage[T1]{fontenc}
\usepackage[dvipsnames]{xcolor}
\usepackage{lmodern}
\usepackage[export]{adjustbox}
\usepackage{subcaption}

\usepackage{amsmath}

\usepackage{stackengine}
\usepackage{mathtools}
\usepackage{stackrel}

\newcommand{\kernel}[2]{ \stackrel[#1]{#2}{  \otimes  }}

\usepackage{amsfonts}
\usepackage{hyperref}
\usepackage{makecell}
\usepackage{amssymb}
\usepackage{lipsum}
\usepackage{multicol}
\usepackage[font=small]{caption}
\usepackage{subcaption}
\usepackage{float}
\usepackage[separate-uncertainty=true]{siunitx}
\usepackage{graphicx}
\usepackage{xfrac}
\usepackage{adjustbox}
\usepackage{standalone}
\usepackage{tikz}
\usepackage{comment}
\usepackage{changepage}
\definecolor{Red}{rgb}{1,0,0}
\definecolor{Green}{rgb}{0,.6,0}
\definecolor{Blue}{rgb}{0,0,1}

\newcommand{\meta}{\mathrm{m}}
\newcommand{\cat}{\mathrm{q}}
\newcommand{\integer}{\mathrm{z}}
\newcommand{\unordered}{\mathrm{q}_\mathrm{u}}
\newcommand{\ordered}{\mathrm{q}_\mathrm{o}}
\newcommand{\standard}{\mathrm{s}}
\newcommand{\continuous}{\mathrm{c}}
\newcommand{\aux}{\mathrm{aux}}

\author{
  \href{https://www.gerad.ca/Charles.Audet/}{Charles Audet}
  \thanks{\url{https://www.gerad.ca/Charles.Audet}}
  \and
  \href{https://www.gerad.ca/fr/people/edward-halle-hannan}{Edward Hall\'e-Hannan}
  \thanks{\href{mailto:edward.halle-hannan@polymtl.ca}{\url{edward.halle-hannan@polymtl.ca}}}
  \and 
  \href{https://www.gerad.ca/Sebastien.Le.Digabel/}{S\'ebastien Le~Digabel}
  \thanks{\url{https://www.gerad.ca/Sebastien.Le.Digabel}}
}

\usepackage{fancyhdr}
\usepackage{csquotes}
\usepackage{xcolor}
\hypersetup{
    colorlinks,
    linkcolor={red!20!black},
    citecolor={blue!20!black},
    urlcolor={blue!80!black}
}

\usepackage[ruled,vlined]{algorithm2e}
\DontPrintSemicolon

\usepackage{amsthm}
\theoremstyle{plain}
\newtheorem{mydef}{Definition}

\begin{document}

% Titre
\title{
{A general mathematical framework for constrained mixed-variable blackbox optimization problems with meta and categorical variables}    
\thanks{GERAD and Department of Mathematics and Industrial Engineering, Polytechnique  Montr\'eal.
This work is supported by the NSERC Alliance grant 544900-19 in collaboration with Huawei-Canada, an IVADO MSc Excellence Scholarship and a Hydro-Qu\'ebec MSc Excellence Scholarship.} }
\maketitle 

%%%%%%%%%%%%%%%%%%%%%
%%%%% Abstract %%%%%%
%%%%%%%%%%%%%%%%%%%%%
\begin{center}
    \text{\Large 
    \textbf{Abstract}}
\end{center}

\begin{adjustwidth}{60pt}{60pt}
A mathematical framework for modelling constrained mixed-variable optimization problems is presented in a blackbox optimization context. 
The framework introduces a new notation and allows solution strategies.
The notation framework allows meta and categorical variables to be explicitly and efficiently modelled, which facilitates the solution of such problems. The new term meta variables is used to describe variables that influence which variables are acting or nonacting: meta variables may affect the number of variables and constraints. The flexibility of the solution strategies supports the main blackbox mixed-variable optimization approaches: direct search methods and surrogate-based methods (Bayesian optimization). The notation system and solution strategies are illustrated through an example of a hyperparameter optimization problem from the machine learning community. \\

\noindent \textbf{Keywords.} Blackbox optimization, derivative-free optimization, mixed-variable optimization, categorical variables, meta variables.
\end{adjustwidth}

%----------------------------------------%
\section{Introduction}
\label{sec:intro}
%----------------------------------------%

This work considers a general constrained optimization problem
\begin{equation}
    \min_{x \in \Omega \subseteq \mathcal{X} } f(x), 
    \label{eq:formulation}
\end{equation}
\noindent where $x \in \mathcal{X}$ is a point that resides in the domain $\mathcal{X}$, $f:\mathcal{X} \to \mathbb{R}$ is the objective function and $\Omega \subseteq \mathcal{X}$ is the feasible set defined by the constraints of the problem.

%----------------------------------------%
\subsection{Context and motivation}
%----------------------------------------%

In blackbox optimization (BBO), the objective and constraint functions are assumed to be blackboxes. In~\cite{AuHa2017}, a mathematical blackbox is defined as: ``\textit{any process that when provided an input, returns an output, but the inner working of the process are not analytically available}''. In general, these blackboxes are computer programs with expensive processes. Hence, in BBO, the objective  and  constraint functions can only provide information through evaluations, which are generally costly. For instance, the only information that the objective function $f$  can provide is the mapping of a point $x\in \mathcal{X}$ to its image $f(x) \in \mathbb{R}$ through a given process, which is typically a computer program. Consequently, derivatives are often inaccessible or too costly to compute. 
Thus, in general, traditional optimization methods cannot be applied to blackbox functions~\cite{AuHa2017}.  

%\footnote{In some cases, $f$ is not a function, but rather a relation. Generally, this subtlety emerges from stochastic routines within the blackbox, thus leading to two different values $f(x)$ for a common point $x$. In this study, the term function is used even for these special cases, which is a slight abuse of terminology. In Section~\ref{sec:example}, an example is described where the objective function $f$ is a performance metric of a multilayer perceptron on a test set. This objective function $f$ implicitly contains  stochastic subroutines in the backpropagation algorithm (training) and the parameters of the model are generally initialized randomly.}

Mixed-variable problems are notoriously hard to tackle in the blackbox optimization community.  This can be partly explained by the presence of meta and categorical variables. Meta variables are a special type of variables that may affect the dimension, the number of constraints and determine which variables are included or excluded in the optimization process. 
These special variables are a cornerstone of this work and are thoroughly defined in Section~\ref{sec:notation}. Moreover, categorical variables are fundamentally difficult to treat since they belong to discrete sets that do not contain any intrinsic metric of distance between the elements and they cannot be relaxed easily in comparison to integer variables. An example of categorical variable is the blood-type of a given person
$x \in \{\text{O-, O+, A-, A, }\ldots \}$
In conjunction, the meta and categorical variables give rise to a substantial challenge in a context of blackbox optimization. In addition, Problem~\eqref{eq:formulation} may contain continuous or integer variables. 

The compact and general formulation of Problem~\eqref{eq:formulation} does not explicitly model mixed-variable problems.
Hence, in order to efficiently tackle these problems, the formulation must be further detailed with a focus on treating the meta variables and the categorical variables. 
A core aspect of this work is to formally define the domain $\mathcal{X}$, which has many implications in the mathematical framework that consists of a notation system and solution strategies. 
The notation framework rigorously models constrained mixed-variable problems in an efficient and unambiguous manner, as well as shines the light on some algorithmic subtleties in the solution of these problems. 
The present work also formalizes solution strategies present in the literature and tackles these problems by being fully compatible with the main blackbox optimization approaches: direct search and surrogate-based Bayesian optimization approaches. The present work does not present any computational experiments,as it focuses on the presentation of the framework, in the same way that the well-known surrogate management framework~\cite{BoDeFrSeToTr99a}, was proposed without experiments. 

%----------------------------------------%
\subsection{Literature review}
\label{sec:literature_review}
%----------------------------------------%

A first framework to treat mixed-variable optimization problems in a context of blackbox optimization is detailed in~\cite{AuDe01a}. The methodology is based on the general pattern search algorithm (GPS) and the variables are partitioned into two components:  discrete and  continuous. The discrete component contains both the quantitative and the qualitative discrete variables, {\em i.e.}, integer variables in $\mathbb{Z}$ as well as  categorical variables. The continuous component contains the continuous variables. Two main ideas emerged from this article. 
First, the continuous space, in which classical continuous blackbox optimization methods can be applied, are generated after fixing the discrete component. Thus, for a fixed discrete component, a continuous space is generated and explored. 
Second, the exploration of the discrete variables space is being done by defining a set of neighbors function $\mathcal{N}$, which is an additional structure to the domain $\mathcal{X}$, such that $\mathcal{N}(x)$ is a set of neighbors of $x \in \mathcal{X}$. With this additional structure a local minimizer $x_{\star}$ is defined so that $x_{\star}$ minimizes the objective function $f$ with respect to the set of neighbors (discrete part) and the continuous space. From the contributions of~\cite{AuDe01a}, a practical application of a thermal insulation optimization problem is treated and optimized~\cite{KoAuDe01a}. In~\cite{AbAuDe2007a}, the filter method is added to the methodology proposed in~\cite{AuDe01a, KoAuDe01a}. This addition enables the methodology to treat general nonlinear constraints. In~\cite{AACW09a}, the methodology based on GPS in~\cite{AbAuDe2007a, AuDe01a, KoAuDe01a} is extended to the mesh adaptive direct search (MADS) algorithm~\cite{AuDe2006}. A rigorous convergence analysis based on~\cite{AbAuDe2007a} was improved by using the Clarke generalized derivatives on the continuous space. In~\cite{AuLeDTr2018}, the MADS algorithm is equipped with a granular mesh called GMesh, which allows the discretization of granular and continuous variables simultaneously. Granular variables are quantitative variables with controlled number of decimals. In particular, GMesh enables to treat integer-continuous problems with the MADS algorithm since integer variables are a special type of granular variables without decimals. 

An important contribution from~\cite{LuPi04a, LuPiSc05a} is the introduction of dimensional variables. These variables affect the number of variables, the number of constraints and the structure of the optimization problem. A point $x$ is partitioned into three components: a dimensional component, a discrete component and a continuous component. The discrete set, where the discrete component belongs, is generated from a fixed dimensional component. Additionally, the continuous space is generated from both a fixed dimensional component and a fixed discrete component. From the partition of an point $x$, a domain and a feasible set are implicitly presented in the formulation of an optimization problem. The present work importantly relies on the contributions from~\cite{LuPi04a, LuPiSc05a}.  

A categorical kernel function is defined in~\cite{Na2021} with the aim of tackling mixed-variable optimization problems with a surrogate approach based on radial basis functions (RBF). The categorical kernel function measures the number of disagreement between two categorical components, where a disagreement is counted when a specific variable of the two compared components is not the same. The surrogate is built upon a composed kernel such that the RBF, centered at some interpolation points, are shifted by the number of categorical disagreements between the fitted point and the interpolation points. The criterion to determine which point is evaluated by the objective function is based on~\cite{RGRegis_CAShoemaker_2007c}. In essence, the criterion has a high value for points that are distant from the previous evaluations (exploration) or points that have promising surrogate-value (intensification). 

Bayesian optimization (BO) has undergone significant development with the recent advent of machine learning. Nowadays, the emerging scientific literature is mainly related to BO based on Gaussian processes (GP), which serves as probabilistic
distribution surrogates~\cite{RaWi06}. The success of BO is explained by an acquisition function that selects which candidate point is to be evaluated. The acquisition function defines a less costly optimization problem with the surrogate. For continuous problems, a well documented acquisition function is the expected improvement ($EI$)~\cite{JoScWe1998}, which provides candidate points in unexplored regions (exploration) and candidate points in promising regions (intensification):  algorithms that applies an $EI$ function on a GP are often referred to as efficient global optimization (EGO) algorithms~\cite{JoScWe1998}. Historically, BO based on GPs was used to tackle continuous blackbox optimization problems. Hence, in practice the integer and categorical variables (one-hot encoded) are often relaxed as continuous variables and rounded afterwards~\cite{GaHe2020}. This naive approach, used in some modern blackbox solvers, often leads to failure such as a mismatch between the points provided by an acquisition function and where the true evaluation takes place, as well as reevaluating some points~\cite{GaHe2020}. Moreover, an important number of additional variables may be generated by the one-hot encoding of categorical variables. In~\cite{RoPaDeClPeGiWy2020}, continuous-categorical optimization problems are modelled with GPs, where a GP surrogate is characterized by kernel composed of tensor products and additions of one-dimensional kernels: an one-dimensional kernel per variable. The one-dimensional kernel of a given categorical variable $x_j \in \{1,2,\ldots, C \}$ is a $C \times C$ matrix, where an element of the matrix is a correlation measure between two categories (classes) of $x_j$. The matrix-kernels for the ordinal and nominal categorical variables are distinguished.
In~\cite{PeBrBaTaGu2021}, the BO framework is extended to tackle mixed-variable optimization problems with continuous, discrete (categorical and integer) and dimensional variables, such as defined in~\cite{LuPi04a, LuPiSc05a}. Again, the GP surrogate is  characterized by a composed kernel built upon products and additions of one-dimensional kernels, each specified by the type of its corresponding variable. Moreover, two approaches are proposed in~\cite{PeBrBaTaGu2021}: multiple surrogates, one surrogate per dimensional component (set of dimensional variables), which separates the main problems into subproblems and a single surrogate with a composed kernel built upon on all variables, including dimensional variables.

In~\cite{MMZ2019}, the authors
combined the user-defined set of neighbors in order to tackle categorical variables in an EGO subproblem. More precisely, a user-defined set of neighbors is randomly defined with a discrete probability distribution based on a GP. Thus, the randomly user-defined set of neighbors serves as a randomized categorical exploration strategy for the EGO subproblem.

% The authors in~\cite{ZhTaChAp2020} argued that most blackbox optimization methods for quantitative-qualitative problems are based on a multiresponse GP with a different response for each combination of levels qualitative factors (each categorical component). The multiresponse essentially comes from the categorical kernel, which is generally built upon tensor products, tensor additions and simplifications in the covariance function.

Covariance functions (kernels) are fundamentally difficult to defined on categorical sets since the distance between two categories (levels) is not defined. To tackle this difficulty, the authors in~\cite{ZhTaChAp2020} proposed to map the categories of each categorical variable to a set of quantitative values that represents some underlying latent unobservable quantitative variable. More precisely, the categories of each categorical variable are mapped to a 2D continuous space: for a given categorical variable, the categories are compared into a 2D space. The quantitative values in the vectors does not have any intrinsic meaning. However, the distance between the values encapsulates some information, since the categories are mapped among themselves in a correlated manner. Mathematically, the mapping is done via a maximum likelihood estimation (MLE) procedure that fits the best multivariate Gaussian distribution of some data. A GP model is then constructed on continuous variables and latent variables. Furthermore, the authors in~\cite{CuLeRoPeDuGl2021} formalized a pre-image problem with a constraint that recovers a categorical component from a vector of continuous latent variables. More technically, a continuous EGO problem is formulated as an augmented Lagrangian with a retrieving constraint on the continuous latent variables. 

The document is organized as follows. First, an example of a mixed-variable optimization problem, taken from the machine learning community, is described in Section~\ref{sec:example}. The example is used throughout the paper to facilitate understanding. Second, the notation system is exhaustively detailed in Section~\ref{sec:notation}. The notation partitions variables in different types, classifies constraint functions, and formally presents their domain and the feasible set. 
Finally, solution strategies are presented in Section~\ref{sec:algorithmic_framework} from the framework perspective.

%----------------------------------------%
\section{Hyperparameter multilayer perceptron example} 
\label{sec:example}
%----------------------------------------%

In order to illustrate the mathematical framework, a simplified constrained hyperparameter optimization problem on a multilayer perceptron (MLP) is detailed throughout the document. Some important hyperparameters are internationally left out, such as the mini-batch size or dropout. The goal of the detailed problem is to model a simple constrained mixed-variable optimization problem in a deep learning context. The objective function is composed of the training and testing of a deep neural network model on a given task. The goal is to find the set of hyperparameters that maximizes a performance score, which is usually a precision score of accuracy on a untested data set. %This problem is compatible with Problem~\eqref{eq:formulation}.

In the example, the MLP is defined to perform regression for inputs with $p \in\mathbb{N}$ continuous features. In other words, the MLP approximates a nonlinear function $h:\mathbb{R}^{p} \to \mathbb{R}$. In order to respect the dimensions of the domain and the codomain of the function $h$, the architecture of the network must have $u_{in}=p$ units in the input layer and $u_{out}=1$ unit in the output layer. The hyperparameters of the MLP are described in Table~\ref{tab:hyperparameter}.

\begin{table}[htb!]
\small
\centering
\begin{tabular}{@{}llcl@{}}
\toprule
\multicolumn{2}{l}{Hyperparameter}               & Variable     & Scope                                        \\ \midrule
\multicolumn{2}{l}{Learning rate}                & $r$          & $]0, 1[$                                 \\
\multicolumn{2}{l}{Activation function}          & $a$          & $ \{$ReLU, Sigmoid$\} $                  \\
\multicolumn{2}{l}{\# of hidden layers}          & $l$          & $\{ 0,1, \ldots, l^{\max} \} $                  \\
\multicolumn{2}{l}{\hspace{0.2cm} \# of units hidden layer $i$} & $u_{i}$      & $ \{u_{i}^{\min} ,u_{i}^{\min}+1, \ldots, u_{i}^{\max} \} $ \\
\multicolumn{2}{l}{Optimizer}                    & $o$          & $ \{$Adam, ASGD$\} $                  \\
\multicolumn{2}{l}{\hspace{0.2cm} if $o=$ ASGD}  &              &                                              \\
       & \hspace{0.4cm} decay                    & $\lambda$    & $]0, 1[$                                     \\
       & \hspace{0.4cm} power update             & $\alpha$     & $]0, 1[$                                     \\
       & \hspace{0.4cm} averaging start          & $t_0$        & $] 1\text{E}3, 1\text{E}8[$                  \\
\multicolumn{2}{l}{\hspace{0.2cm} if $o=$ Adam} &              &                                              \\
       & \hspace{0.4cm} running average 1        & $\beta_1$    & $]0, 1[$                                     \\
       & \hspace{0.4cm} running average 2        & $\beta_2$    & $]0, 1[$                                     \\
       & \hspace{0.4cm} numerical stability      & $ \epsilon $ & $]0, 1[$ 
       \\ \bottomrule
\end{tabular}
\caption{Hyperparameters of the MLP.}
\label{tab:hyperparameter}
\end{table}

The index $i$ in $u_i$ represents the $i$-th hidden layer. The number of units in the hidden layers are grouped in the vector $u(l)=(u_1,u_2,\ldots,u_l)$, where $l$ is the number of hidden layers. The situation where there are no hidden layer is modeled by setting $l=0$. In that case, the variables $u_i$ are said to be nonacting, which signifies that the variables $u_i$ are not part of the optimization problem when $l=0$. The terminology of acting and nonacting are further detailed in Section~\ref{sec:meta_component} and Section~\ref{sec:roles}.

Depending on the choice of the optimizer, different hyperparameters are involved. Indeed, in Table~\ref{tab:hyperparameter} the optimizers do not share the same continuous hyperparameters. A given optimizer leads to different variables in the problem. For example, the variable decay $\lambda$ is only part of the problem (acting) if $o=\text{ASGD}$. This consideration is important and will be discussed throughout the document, but notably in Section~\ref{sec:meta_component} that focuses on meta variables.  

\begin{figure}[htb!]
\centering
  \includegraphics[width=0.60 \linewidth]{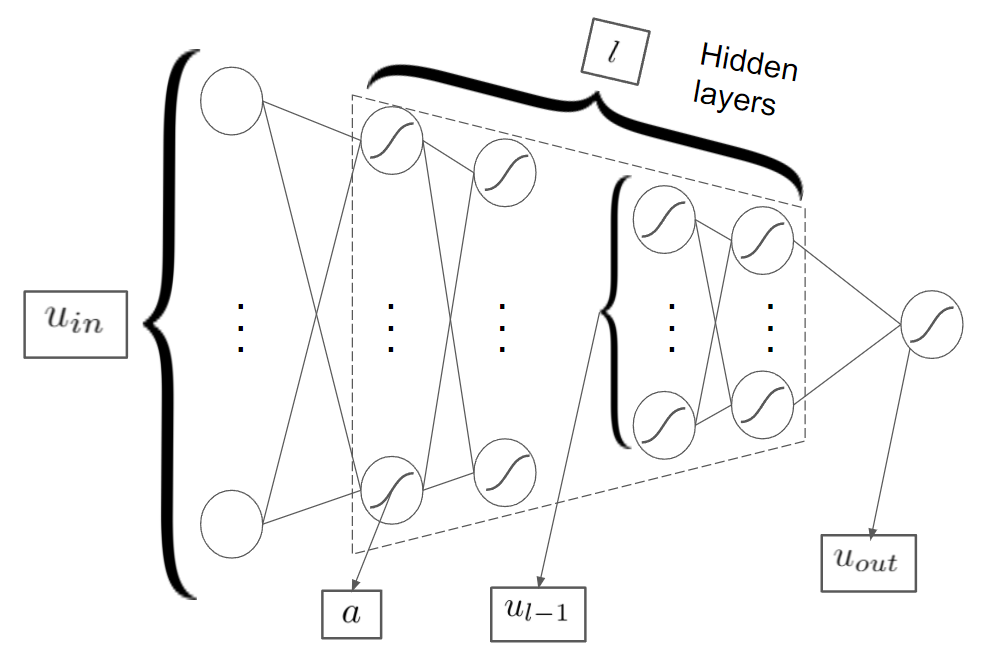}
  \caption{MLP of the hyperparameter problem (see Table~\ref{tab:hyperparameter}).}
  \label{fig:MLP}
\end{figure}

The first constraint of the problem imposes that the sum of the units in all the hidden layers does not exceed an upper bound $\hat{u} \in \mathbb{N}$, such that $ \sum_{i=1}^{l} u_i \leq \hat{u}$. The other constraints are
$ u_{i} \leq u_{i-1} \ \forall i \in \{2, 3,\ldots, l \}$ and they impose that the number of units in subsequent hidden layers are less than or equal, which may help reduce the number of units. These artificial constraints are imposed to illustrate the notation (see Section~\ref{sec:worked_example_notation}). An example of a possible architecture of the MLP is schematized in Figure~\ref{fig:MLP}.

%----------------------------------------%
\section{Notation framework}
\label{sec:notation}
%----------------------------------------%

This section contains the fundamental mathematical definitions that allow modelling mixed-variable problems. In Section~\ref{sec:variable}, the mathematical objects that define the variables (point and components) are described. Subsequently, the domain $\mathcal{X}$ is detailed in Section~\ref{sec:domain}. Then, the feasible set $\Omega  \subseteq \mathcal{X}$ is precised in Section~\ref{sec:set_feasible}. Finally, the content in Sections~\ref{sec:variable} to~\ref{sec:set_feasible} is discussed within the MLP example in Section~\ref{sec:worked_example_notation}.  

%-----------------------------------------------%
\subsection{Variables and components of a point}
\label{sec:variable}
%-----------------------------------------------%

The goal of an optimization algorithm is to find a feasible point $x_{\star}$ that minimizes the objective function $f$. In a mixed-variable optimization context, it is necessary to formally define how a point is partitioned into different components.

\begin{mydef}[Components of a point]
A point $x=(x^{\meta},x^{\cat},x^{\standard})$ is partitioned into three components:
\begin{itemize}
    \item a meta component $x^{\meta}$;
    
    \item   a categorical component $x^{\cat}=(x^{\unordered},x^{\ordered})$, which itself is partitioned into the unordered categorical (nominal) component $x^{\unordered}$ and the ordered categorical (ordinal) component $x^{\ordered}$;
    
    \item   a standard component $x^{\standard}=(x^{\integer},x^{\continuous})$, which itself is the fusion of the integer component $x^{\integer}$ and the continuous component $x^{\continuous}$.
\end{itemize}

\noindent For each $t \in \{\meta, \cat, \standard, \unordered, \ordered, \integer, \continuous  \}$, the component $x^t$ is a vector containing $n^t \in \mathbb{N}$ variables of type $t$:
\begin{equation}
    x^t = (x^t_1, x^t_2, \ldots, x^t_{n^t}).
    \label{eq:component}
\end{equation}
\label{def:point}
\end{mydef}

The fusion of the integer and continuous components into the standard component $x^{\standard}$ is justified by several reasons. 
In practice these variables are generally optimized with standard methods. Moreover, some blackbox optimization algorithms have the ability to simultaneously optimize integer and continuous variables. Thus, it is convenient to group these variables to lighten the notation. However, the standard component $x^{\standard}=(x^{\integer}, x^{\continuous})$ can easily be partitioned into its two components if necessary.

The meta, standard and categorical components, as well as their corresponding variables, are respectively discussed in Sections~\ref{sec:meta_component},~\ref{sec:categorical_component}
and~\ref{sec:standard_component}. Additionally, the motivations behind the compact partition $x=(x^{\meta},x^{\cat},x^{\standard})$ and the complete partition $x=(x^{\meta},x^{\unordered},x^{\ordered}, x^{\integer},x^{\continuous})$ are discussed and illustrated in Section~\ref{sec:variable_types}. Finally, in Section~\ref{sec:roles}, the roles of variables and constraints are introduced in order to define more clearly the domain $\mathcal{X}$ in Section~\ref{sec:domain}.

%--------------------------------------------------%
\subsubsection{Meta component and decree property}
\label{sec:meta_component}
%--------------------------------------------------%

The meta component $x^{\meta}$ contains variables having the decree property, which are called meta variables. The decree property is a special property that only meta variables possess. 
The property determines if some variables or constraints are either acting or nonacting. 
The term acting indicates that the variable or constraint is part of the problem and, on the contrary, the term nonacting indicates that the variable or constraint is not part of the problem.  
More precisely, an acting variable is a decision variable that is included in the domain in which the optimization process is deployed. 
An acting constraint is a constraint function that defines the feasible set that contains feasible solutions. 

The decree propriety is  attributed to the meta component $x^{\meta}$, since it contains the meta variables. 
Concretely, meta variables may affect the number of variables (dimension) or the number of constraints. With the MLP example, the number of hidden layers $l$ affect the dimension and the number of constraints of the problem, whereas the optimizer $o$ does not. Indeed, $l$ affect the number of variables (dimension) since it decrees the units $u_i$ in the hidden layers $i \in \{1,2, \ldots, l^{\max} \}$, such that $u_i \in \{u_1,u_2,\ldots, u_l \}$ are acting variables and $u_i \in \{u_{l+1}, u_{l+2}, \ldots, u_{l^{\max}} \}$ are nonacting. Moreover, $l$ also decrees the corresponding constraints $u_{i} \leq u_{i-1} \ \forall i \in \{2,3,\ldots, l\}$, thus affecting the number of constraints. 
Both optimizers ASGD and Adam from Table~\ref{tab:hyperparameter} decree three continuous hyperparameters, which does not affect the dimension nor the number of constraints. However, the optimizer decrees some variables. For instance, the decay $\lambda$ is only an acting variable if $o=\text{ASGD}$. 

In that regards, meta variables are a generalization of the strictly discrete dimensional variables defined in~\cite{LuPi04a, LuPiSc05a}. First of all, meta variables do not necessarily affect the dimension, in comparison to dimensional variables. Secondly, meta variables can be of any type. For example, a problem could contain a continuous variable frequency that takes its value in the visible spectrum (continuous scope). The visible spectrum could be partitioned into the three intervals that represent the red-blue-green colors. Finally, the frequency could decree some variables or constraints depending in which interval (color) it belongs to. In that particular example, the frequency is a meta-continuous variable.

Additionally, the terminology {\em dimensional} used in~\cite{LuPi04a, LuPiSc05a, PeBrBaTaGu2021} is avoided, since it is used in physical sciences and engineering to describe quantities such as the velocity, mass and time. Many of blackbox mixed-variable optimization problems come from these disciplines.

%--------------------------------------------------%
\subsubsection{Roles of variables and constraints}
\label{sec:roles}
%--------------------------------------------------%

The present section introduces {\em roles} of variables and constraints.
They are introduced for two reasons : 1) they facilitate the comprehension of the influence of meta variables on the other variables of type $t\in \{ \cat, \standard, \unordered, \ordered, \integer, \continuous \}$ and constraints; 2) the definition of the domain $\mathcal{X}$ in Section~\ref{sec:domain} and the feasible set $\Omega$ in Section~\ref{sec:set_feasible} are made clearer. 
In essence, the roles of variables and constraints consist of additional terminologies that help elucidate some subtleties of the mathematical framework.

The role of a variable must not  be confused with its variable type. 
In addition to its type $t\in \{\meta,  \cat, \standard, \unordered, \ordered, \integer, \continuous \}$, each variable takes a single role amongst meta, decreed or global. 
A constraint takes a single role amongst decreed or global. 

The roles of meta variables is simply their meta type : meta is both a variable type and a role. 
The role of meta variables is to decree their decreed variables or constraints. More precisely, some variables or constraints may be acting or nonacting accordingly to some specific meta variables. These variables or constraints are said to be decreed. They are called decreed variables and decreed constraints. In the MLP example of Table~\ref{tab:hyperparameter}, the decay $\lambda$ is an acting variable if the optimizer $o=\text{ASGD}$ and is nonacting otherwise. Thus, the optimizer $o$ is a meta variable and the decay $\lambda$ is a decreed variable. The decay $\lambda$ is also said to be decreed by the optimizer $o$. Conceptually, the role of decreed variables or constraints is to be acting or nonacting with respect to to their specific meta variables\footnote{Note that a decreed variable can not be a meta variable: this modeling choice is being done to allow only one instance of meta variables. In other words, meta variables that decrees meta variables are not allowed, since it is very uncommon to encounter such problems in practice and it would complexify the notation even more.}.

Additionally, the decree property has been attributed to meta component $x^{\meta}$ in Section~\ref{sec:meta_component}. Thus, the meta component $x^{\meta}$ has an important role in which it decrees all the decreed variables, since the meta component $x^{\meta}$ contains all the meta variables.

The last role is the simplest one. Global variables or constraints are always acting and do not possess the decree property. In other words, the global variables or constraints are not meta and are always part of the problem. Conceptually, global variables or constraints have an empty role in the sense that they do not influence and are not influenced by other variables or constraints. In the MLP example of Table~\ref{tab:hyperparameter}, the activation function $a$ is a global variable since it is not decreed by any variable and it does not decree other variable. 

\begin{figure}[htb!]
\centering
  \includegraphics[width=0.75 \linewidth]{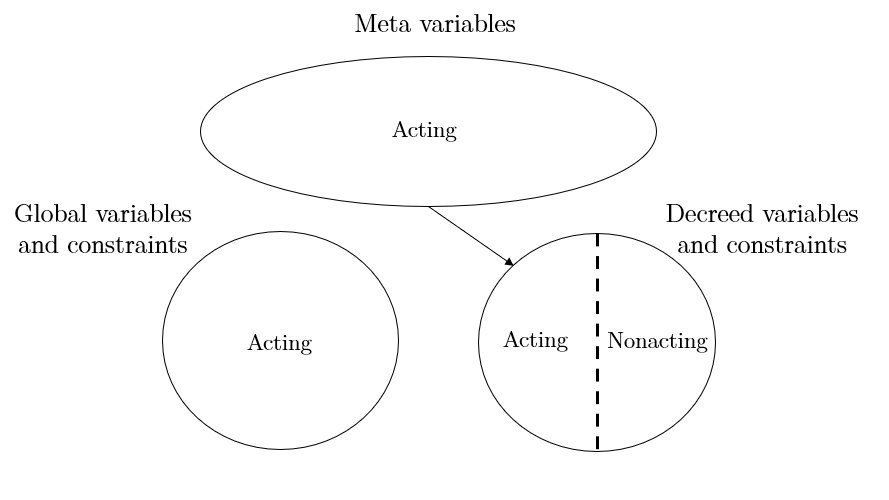}
  \caption{Role classification of variables and constraints.}
  \label{fig:roles}
\end{figure}

Figure~\ref{fig:roles} summarizes the roles of variables and constraints by illustrating some important concepts. First, the arrow symbolizes that meta variables decree some variables or constraints. Second, nonacting variables and constraints are a subset of decreed variables and constraints. This implies that nonacting variables or constraints are necessarily decreed. Third, the global variables and constraints are disjoint, which indicates that they are unaffected by the meta variables.

%--------------------------------------------------%
\subsubsection{Categorical component}
\label{sec:categorical_component}
%--------------------------------------------------%

The categorical component $x^{\cat}$ contains qualitative variables, known as categorical variables, that are not meta: a categorical variable may be decreed or global. Categorical variables are discrete variables that take qualitative values called categories. More precisely, a categorical variable $x^{\cat}_j$ has $c_j$ categories, such that  $x^{\cat}_j \in \{1,2,\ldots, c_j \}$. 

Categorical variables can be unordered or ordered. Unordered categorical variables are known as nominal variables ({\em e.g.}, the blood-type) and they are contained in the unordered categorical component $x^{\unordered}$. Note that binary variables are nominal variables. Subsequently, ordered categorical variables are also known as ordinal variables and they are contained in the  categorical ordered component $x^{\ordered}$. The size of a pizza $x \in \{\text{small},\text{ medium}, \text{ large} \}$ is an ordinal variable, since the categories are ordered from small to large. Although the ordinal variables belong to ordered sets, distances between the ordinal variables are inherently unknown: ``[\ldots] \textit{there is an ordering between the values, but no metric notion is appropriate}''~\cite{hastie01statisticallearning}. For short, the terms nominal and ordinal are prioritized over unordered categorical and ordered categorical respectively.  

The categorical component $x^{\cat}=(x^{\unordered}, x^{\ordered})$ is composed of the nominal component $x^{\unordered}$ and the the ordinal component $x^{\ordered}$, which respectively contains the nominal variables and the ordinal variables that are not meta. In some cases, it might be beneficial to exploit the order of an ordinal variable, motivating the partition of the categorical component into nominal and ordinal components. For instance,~\cite{RoPaDeClPeGiWy2020} used different kernels for ordinal and nominal components. Moreover, a direct search exploration strategy could be generically implemented with a previous and next element mechanism for an ordinal set.

In previous works~\cite{AbAuDe2007a, AuDe01a, KoAuDe01a}, meta variables were included in the categorical variables; it is an important distinction from this work.

%--------------------------------------------------%
\subsubsection{Standard component}
\label{sec:standard_component}
%--------------------------------------------------%

The standard $x^{\standard}$ component contains discrete and continuous quantitative variables that are not meta variables: a standard variable may be decreed or global. Formally, the standard component $x^{\standard}$ contains variables that belong to intrinsically ordered sets for which a metric of distance is intuitively definable. Simply put, the standard component $x^{\standard}$ contains the integer variables and the continuous variables.
    
The integer component $x^{\integer}$ exclusively contains discrete quantitative variables, called integer variables, that are not meta. Unlike the categorical variables, integer variables are always ordered and belong to sets with appropriate metric notions~\cite{hastie01statisticallearning}. The decision to separate the discrete variables into the the categorical component and the integer component differs from some of the current literature. Indeed, in~\cite{AACW09a, LuPiSc05a, PeBrBaTaGu2021} the discrete component contains both the categorical and the integer variables. Thus, categorical and integer variables are not clearly distinguished: some useful mathematical properties of the integer variables might not be exploited at their fullest. In that regard, integer programming is a well developed optimization field that exploits the properties of the integer variables. In practice, this strengthens the separation of the integer variables from the categorical ones, since integer programming techniques could be implemented in the algorithmic framework to treat the integers variables. 

The continuous component $x^{\continuous}$ contains continuous variables that are not meta. Continuous variables have many properties that are generally exploited in a context of blackbox optimization.  

%--------------------------------------------------%
\subsubsection{Variable type classification}
\label{sec:variable_types}
%--------------------------------------------------%

\begin{figure}[htb!]
\centering
  \includegraphics[width=1.00 \linewidth]{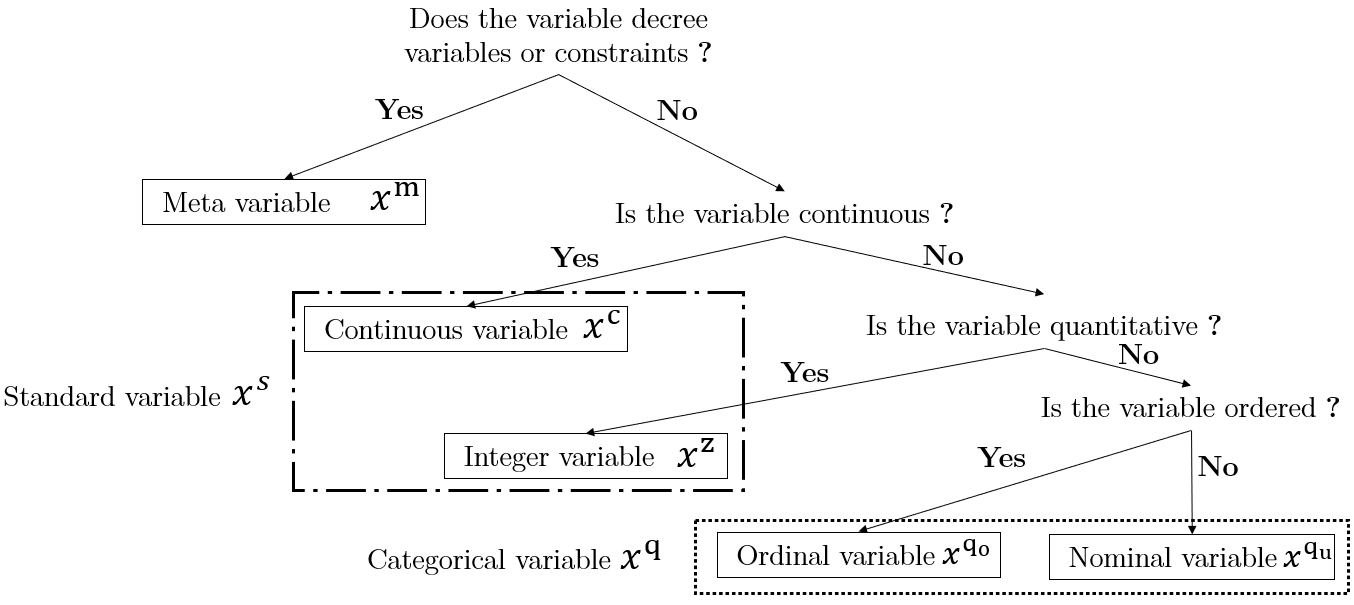}
  \caption{Variable type classification tree chart.}
  \label{fig:variables_id}
\end{figure}

Figure~\ref{fig:variables_id} shows a tree chart that classifies a variable by their type in the proposed mathematical framework. The first question identifies meta variables, the second determines the continuous variables, the third distinguishes integer from categorical variables and the last one separates ordinal from nominal variables. The first question also imply that continuous, integer, ordinal and nominal variables are not meta variables. The doted box in the middle illustrates that standard variables contains the continuous and integer variables, whereas the doted box in the bottom exhibits that two types of categorical variable, which are ordinal and nominal variable.

%Figure~\ref{fig:variables_id_set} is a formal version of Figure~\ref{fig:variables_id}, which both describes the full partition of a point $x~=~(x^{\meta}, x^{\unordered}, x^{\ordered}, x^{\integer}, x^{\continuous})$. More precisely, Figure~\ref{fig:variables_id_set} is a Venn diagram of the set properties of each variable type (contained in their component). Forthwith, the meta variables are separated from the not meta variables with the two large rectangles. Meta variables can be of any type, thereby they may reside anywhere in the right rectangle. The classical set properties, such as order and distance, are organized from the inner circle (continuous) to the outer region contained in the rectangle (unordered). From inner to outer represents from most properties to the least properties. Finally, the standard variables resides in the gridded area and the categorical variables resides in the grey area. 

%\begin{figure}[htb!]
%\centering
%  \includegraphics[width=0.60 \linewidth]{images/variables_properties.png}
%  \caption{Classification of variable types by their set properties and the decree property (meta variables).}
%  \label{fig:variables_id_set}
%\end{figure}

%The full partition thoroughly separates variables into the properties of their respective set.

Mathematically, the partition of a point complete partition $x=(x^{\meta},x^{\unordered},x^{\ordered},x^{\integer},x^{\continuous})$, displayed in Figure~\ref{fig:variables_id}, offers flexibility and extracts most mathematical information accessible to facilitate the optimization process: the modelling choices for the partitions are motivated by these considerations. The compact partition $x=(x^{\meta},x^{\cat},x^{\standard})$ implicitly contains the same information and flexibility of the full partition. However, the compact partition alleviates the notation, which is why it is mostly used throughout this work.

%--------------------------------------------------%
\subsection{Domain} 
\label{sec:domain}
%--------------------------------------------------%

At this stage, the variables have been: 1) classified into different types; 2) organized into components, which forms a partition of a point $x$; 3) attributed roles. The next step is to define the domain $\mathcal{X}$ of the objective function $f:\mathcal{X}\to \mathbb{R}$, such that a point $x \in \mathcal{X}$ resides in that set.     

\begin{mydef}[Domain]
\noindent The domain of objective function is defined by: 
\begin{align}
\begin{split}
    \mathcal{X} =  \big\{ \ (x^{\meta}, x^{\cat}, x^{\standard}) \ \ : \ \
    &x^{\meta} \in \mathcal{X}^{\meta}, \\
    &x^{\cat} \in \mathcal{X}^{\cat}(x^{\meta}), \\
    &x^{\standard} \in \mathcal{X}^{\standard}(x^{\meta}) \ \big\}
\end{split}
\end{align}
\noindent where $\mathcal{X}^{\meta} \subseteq \mathbb{M}^{n^{\meta}} $ is the meta set, $\mathcal{X}^{\cat}(x^{\meta}) \subseteq \mathbb{Z}^{n^{\cat}(x^{\meta})} $ is the parametrized categorical set and $\mathcal{X}^{\standard}(x^{\meta})\subseteq \mathbb{Z}^{n^{\integer}(x^{\meta})}~\times~\mathbb{R}^{n^{\continuous}(x^{\meta})} $ is the parametrized standard set.
\label{def:domain}
\end{mydef}

%Two subtleties regarding Definition~\ref{def:domain} must be explained. First, the trivial bound constraints of each variable are contained in the domain $\mathcal{X}$. Secondly, 
The dependencies of the parametrized categorical set $\mathcal{X}^{\cat}(x^{\meta})$ and parametrized standard set $\mathcal{X}^{\standard}(x^{\meta})$ are defined through a parametrization with respect to the meta component $x^{\meta}$. These parametrizations are a direct consequence of the decree property of the meta component $x^{\meta}$. 

\begin{mydef}[Parametrized set]
A parametrized set $\mathcal{X}^{t}(x^{\meta})$ of type $t\in \{\cat, \standard, \unordered, \ordered, \integer, \continuous \}$ is the set that contains all the components of type $t$, such that
\begin{align}
\begin{split}
    \mathcal{X}^{t}(x^{\meta}) = 
    \Big \{ \ x^t = \left(x^t_1, x^t_2, \ldots, x^t_{n^t(x^{\meta}) }\right) \ : \ x_i^t \in S_i^t \text{ is an acting variable }  \forall i \in I^t(x^{\meta}) \Big \}
\end{split}
\end{align}
where $S_i^t$ is the scope of the acting variable $x_i^t$ and $I^t(x^{\meta})=\{ 1,2,\ldots, n^t(x^{\meta}) \}$ is the set of indices of the acting variables $x_i^t$, which are either global or decreed by the meta component $x^{\meta} \in \mathcal{X}^{\meta}$. 
\label{def:parametrized_set}
\end{mydef}

From Definition~\ref{def:parametrized_set}, it follows that a component $x^{t} \in \mathcal{X}^{t}(x^{\meta})$ contains only the acting variables of type $t$. The nonacting variables are not contained in the component $x^{t} \in \mathcal{X}^t(x^{\meta})$. Recall that nonacting variables are necessarily decreed variables, whereas acting variables may be global or decreed. Hence, in the component $x^t \in \mathcal{X}^t(x^{\meta})$, some acting variables contained may be decreed by the meta component $x^{\meta}\in \mathcal{X}^{\meta}$, which justifies the parametrization of the set $\mathcal{X}^{t}(x^{\meta})$.  

Two additional remarks follow. First, the meta variables are always acting variables, thus the meta set $\mathcal{X}^{\meta}$ has no dependency. Secondly, a parametrized set $\mathcal{X}^t(x^{\meta})$ is a subset of the set that contains all possible components $\mathcal{X}^t$, such that

\begin{equation}
    \mathcal{X}^t(x^{\meta}) \ \subseteq \ \mathcal{X}^t \ = \bigcup_{x^{\meta} \in \mathcal{X}^{\meta}} \mathcal{X}^t(x^{\meta}),
    \label{eq:set}
\end{equation}

\noindent where $t\in \{\cat, \standard, \unordered, \ordered, \integer, \continuous \}$. A component $y^t$ is said to be incompatible with the meta component $x^{\meta}$, if $y^t \in \mathcal{X}^t$ and $y^t \not \in \mathcal{X}^t(x^{\meta})$ (more compactly, $y^t \in \mathcal{X}^t \setminus \mathcal{X}^t(x^{\meta})$). Again, a set $\mathcal{X}^t$ contains all possible components, hence it must contains these incompatibles components. 

In the MLP example, a continuous component $y^{\continuous}$ that contains the decay $\lambda$ (see Table~\ref{tab:hyperparameter}) is incompatible with the meta component $x^{\meta}$. Indeed, if $o=\text{Adam}$ and $y^{\continuous}$ is a continuous component that contains the decay, then $y^{\continuous} \in \mathcal{X}^{\continuous}$ and $y^{\continuous} \not \in \mathcal{X}^{\continuous}(l,\text{Adam})$.

Moreover, if all variables of type $t\in \{\cat, \standard, \unordered, \ordered, \integer, \continuous \}$ are global variables (unaffected by the meta component), then no parametrization is necessary, such that $\mathcal{X}^t(x^{\meta})=\mathcal{X}^t$. 

The subtlety of incompatible component when some variables are decreed is important in Definition~\ref{def:domain} of the domain $\mathcal{X}$. 
It explains why the domain $\mathcal{X}$ from in Definition~\ref{def:domain} is formulated with a categorical parametrized set $\mathcal{X}^{\cat}(x^{\meta})$ and a standard parametrized set $\mathcal{X}^{\standard}(x^{\meta})$ instead of the categorical set $\mathcal{X}^{\cat}$ and a standard set $\mathcal{X}^{\standard}$. Indeed, for a given meta component $x^{\meta} \in \mathcal{X}^{\meta}$, the categorical and standard components reside in their parametrized sets, such that $x^{\cat} \in \mathcal{X}^{\cat}(x^{\meta}$) and $x^{\standard} \in \mathcal{X}^{\standard}(x^{\meta})$, in order to take into account that some categorical or standard variables may be decreed by the given meta component $x^{\meta} \in \mathcal{X}^{\meta}$.

Moreover, the meta component $x^{\meta}$ may affect the dimension $n^t(x^{\meta})$ of the component $x^t \in \mathcal{X}^t(x^{\meta})$. Indeed, some acting variables of type $t$ contained in the component $x^t \in \mathcal{X}^t(x^{\meta})$ may be decreed by the meta component $x^{\meta}$, thus the number of acting variables in this component may vary with the meta component $x^{\meta}$. In simpler terms, the dimension of the component $x^{t}$ may vary with the meta component $x^{\meta}$. Hence, the dimension of the component $x^t$ is a function $n^t: \mathcal{X}^{\meta} \to \mathbb{N}$. Notably in the MLP example in Table~\ref{tab:hyperparameter}, the number of hidden layers $l$ decrees the number of units $u_i$ in the hidden layers, which affects the number of integer variables. Thus, the dimension of the integer component $x^{\integer} \in \mathcal{X}^{\integer}(x^{\meta})$ is determined by the meta component $x^{\meta}$.

%--------------------------------------------------%
\subsubsection{Alternative formulation of the domain}
\label{sec:domain_algo}
%--------------------------------------------------%

The domain $\mathcal{X}$ defined in Definition~\ref{def:domain} offers little insight regarding the visualization and construction of the domain $\mathcal{X}$, especially regarding the parametrized categorical set $\mathcal{X}^{\cat}(x^{\meta})$ and the parametrized standard set $\mathcal{X}^{\standard}(x^{\meta}$). Hence, a more visual and algorithmic formulation of the domain, based on~\cite{AACW09a, AbAuDe2007a}, is proposed: 

\begin{equation}
    \mathcal{X} \ =  \bigcup_{x^{\meta} \in \mathcal{X}^{\meta} } \left( \{ x^{\meta} \} \times  \bigcup_{x^{\cat} \in \mathcal{X}^{\cat}(x^{\meta})} \Big( \{ x^{\cat} \} \times   \mathcal{X}^{\standard}(x^{\meta})  \Big) \right).
    \label{eq:domain2}
\end{equation}

\noindent Following the same logic as in Equation~\eqref{eq:domain2}, the parametrized standard set $\mathcal{X}^{\standard}(x^{\meta})$ is formulated as:

\begin{equation}
    \mathcal{X}^{\standard}(x^{\meta}) \ \  = \ \ \mathcal{X}^{\integer}(x^{\meta}) \times \mathcal{X}^{\continuous}(x^{\meta}) \ \ = \  \bigcup_{x^{\integer} \in \mathcal{X}^{\integer}(x^{\meta})} \Big( \{ x^{\integer} \} \times   \mathcal{X}^{\continuous}(x^{\meta})  \Big)
    \label{eq:met2}
\end{equation}

\begin{figure}[htb!]
\begin{subfigure}{0.62\textwidth}
  \includegraphics[width=\linewidth]{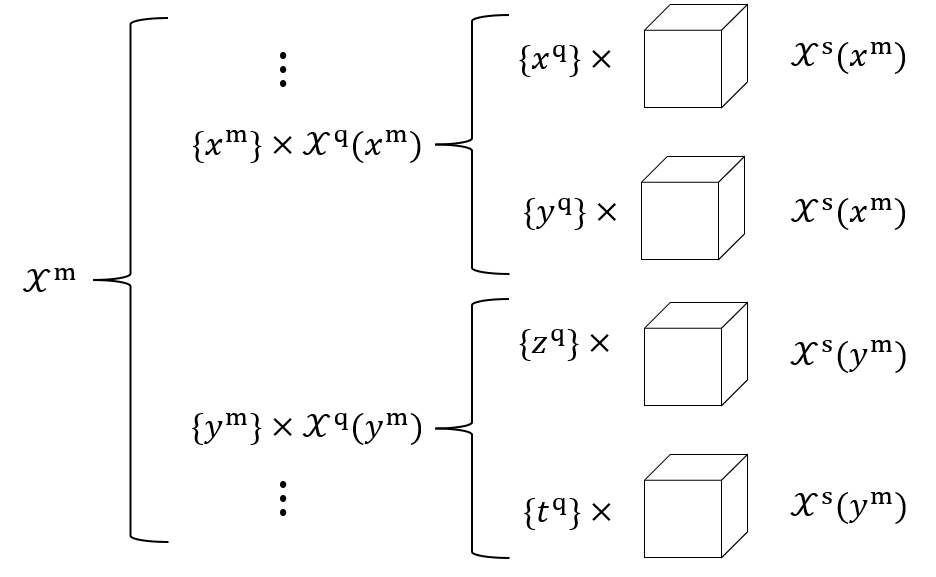}
  \subcaption{Visualization of Equation~\eqref{eq:domain2}.}
  \label{fig:space}
\end{subfigure}
\begin{subfigure}{0.38\textwidth}
  \includegraphics[width=1\linewidth]{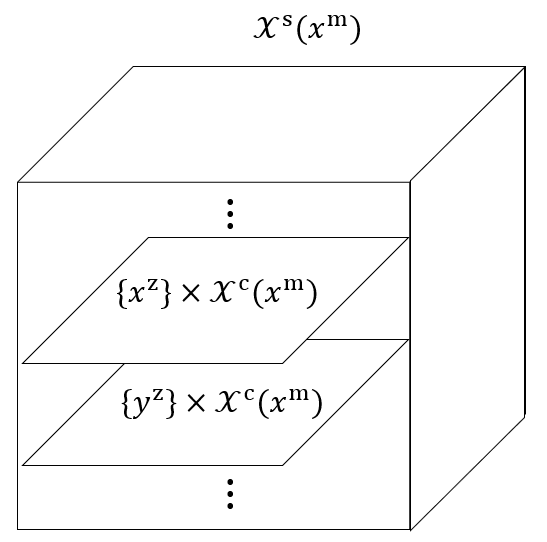}
  \subcaption{Visualization of Equation~\eqref{eq:met2}.}
  \label{fig:box}
\end{subfigure}
\caption{Visualization of the domain $\mathcal{X}$.}
\end{figure}

Schematically, the parametrized standard set $\mathcal{X}^{\standard}(x^{\meta})$ can be visualized as the union of multiple layers, where each layer is a Cartesian product of a parametrized continuous set $\mathcal{X}^{\continuous}(x^{\meta})$ with an integer component $x^{\integer} \in \mathcal{X}^{\integer}(x^{\meta})$, which is illustrated in Figure~\ref{fig:box}. In Figure~\ref{fig:box}, each layer shares the same continuous set $\mathcal{X}^{\continuous}(x^{\meta})$, whereas each layer has a distinct integer component $x^{\integer} \in \mathcal{X}^{\integer}(x^{\meta})$. The standard set $\mathcal{X}^{\standard}(x^{\meta})$ is represented as a box containing all the possible unions described in Equation~\eqref{eq:met2}.  

A visualization of the entire domain $\mathcal{X}$ can be built upon the abstraction of the standard set $\mathcal{X}^{\standard}(x^{\meta})$ illustrated in Figure~\ref{fig:box}. In Figure~\ref{fig:space}, the standard sets are represented as small boxes, following the abstraction from Figure~\ref{fig:box}. The left-curly brackets represents the unions in the Equation~\eqref{eq:domain2}, from left to right. Furthermore, the formulation of the domain $\mathcal{X}$ in Equation~\eqref{eq:domain2} may be understood and visualize as an explicit enumeration of all the possible points, similarly to a set of all possible components $\mathcal{X}^t$ in Equation~\eqref{eq:set}.

The domain $\mathcal{X}$, defined in Definition~\ref{def:domain} or alternatively formulated in Equation~\eqref{eq:domain2}, is composed of the meta set $\mathcal{X}^{\meta}$, the parametrized categorical set $\mathcal{X}^{\cat}(x^{\meta})$ and the parametrized standard set $\mathcal{X}^{\standard}(x^{\meta})$. Hence, few details about these important sets are given in the following sections.

%--------------------------------------------------%
\subsubsection{Meta set}
%--------------------------------------------------%

The number of variables in the meta component is denoted $n^{\meta} \in \mathbb{N}$. The meta component $x^{\meta}$ belongs to the meta set $\mathcal{X}^{\meta} \subseteq \mathbb{M}^{n^{\meta}}$, which contain all the meta component $x^{\meta}$. In comparison to a parametrized set, the meta set $\mathcal{X}^{\meta}$ is static, since the meta variables are always acting variables. This also implies that the meta component $x^{\meta}$ has a fixed dimension $n^{\meta} \in \mathbb{N}$. Moreover, the set $\mathbb{M}^{n^{\meta}}$ is a mixed set consisting of Cartesian products, such that
\begin{equation}
\mathbb{M}^{n^{\meta}} = \mathbb{Z}^{n^{\meta_{\cat}}} \times \mathbb{Z}^{n^{\meta_{\integer}}} \times \mathbb{R}^{n^{\meta_{\continuous}}},
\end{equation}
\noindent where $n^{\meta}=n^{\meta_{\cat}} + n^{\meta_{\integer}} + n^{\meta_{\continuous}} $ is the number of meta variables, $n^{\meta_{\cat}}$ is the number of meta-categorical variables, $n^{\meta_{\integer}}$ is the number of meta-integer variables and $n^{\meta_{\continuous}}$ is the number of meta-continuous variables. In particular, note that $\mathbb{M}^{n^{\meta}}~=~\mathbb{Z}^{n^{\meta_{\cat}}} \times \mathbb{Z}^{n^{\meta_{\integer}}}~=~\mathbb{Z}^{n^{\meta}}$ in the case where meta variables are strictly dimensional variables (discrete variables)~\cite{LuPi04a, LuPiSc05a}. This case is common in practice. 
%--------------------------------------------------%
\subsubsection{Parametrized categorical set}
%--------------------------------------------------%

The categories of each categorical variable can be mapped with a bijection to a subset of $\mathbb{Z}$. Hence, without any loss of generality, the parametrized categorical set $\mathcal{X}^{\cat}(x^{\meta})$ is considered to be a subset of $\mathbb{Z}^{n^{\cat}(x^{\meta})}$. However, this bijection does not imply that a metric notion is appropriate~\cite{hastie01statisticallearning}. In other words, this bijection is only useful in terms of algorithmic implementations.

Since the categorical variable $x_j^{\cat}$ takes values from the set
    $\{1,2,\ldots, c_j \}$, the parametrized categorical set $\mathcal{X}^{\cat}(x^{\meta})$ is defined as
\begin{equation}
     \mathcal{X}^{\cat}(x^{\meta}) = \prod_{j=1}^{n^{\cat}(x^{\meta})} \{1,2,\ldots, c_j \}.
\end{equation}
It may also be expressed as the Cartesian product between the parametrized unordered and ordered sets
\begin{equation}
    \mathcal{X}^{\cat}(x^{\meta}) = \mathcal{X}^{\unordered}(x^{\meta}) \times \mathcal{X}^{\ordered}(x^{\meta}) =  \prod_{i=1}^{n^{\unordered}(x^{\meta})} \{1,2,\ldots, c_i \} \times \prod_{i=j}^{n^{\ordered}(x^{\meta})} \{1,2,\ldots, c_j \}   , 
\end{equation}

\noindent which outlines the distinction between nominal and ordinal variables. 

%--------------------------------------------------%
\subsubsection{Parametrized standard set}
%--------------------------------------------------%

The parametrized standard set $\mathcal{X}^{\standard}(x^{\meta})$ is a compact notation that describes a direct Cartesian product of the parametrized integer and continuous sets:
\begin{equation}
\mathcal{X}^{\standard}(x^{\meta}) \ = \  \mathcal{X}^{\integer}(x^{\meta}) \times \mathcal{X}^{\continuous}(x^{\meta}) \subseteq \ \mathbb{Z}^{n^{\integer}(x^{\meta}) } \times \mathbb{R}^{n^{\continuous}(x^{\meta})} ,
\end{equation}
\noindent where $\mathcal{X}^{\integer}(x^{\meta}) \subseteq \mathbb{Z}^{n^{\integer}(x^{\meta}) }$ is the parametrized integer set and $\mathcal{X}^{\continuous}(x^{\meta}) \subseteq \mathbb{R}^{n^{\continuous}(x^{\meta})} $ is the parametrized continuous set. 

Again, the compact notation for the parametrized standard set $\mathcal{X}^{\standard}(x^{\meta})$ is particularly interesting when for algorithms that optimize simultaneously the integer and continuous variables is employed.

%--------------------------------------------------%
\subsection{Feasible set}
\label{sec:set_feasible}
%--------------------------------------------------%

The constraints are separated into two roles, the global and decreed constraints. 
Global constraints are always acting, whereas decreed constraints may be acting or nonacting depending on the meta variables. 
The decreed constraints lead to the following definition.

\begin{mydef}[Set of decreed acting constraints]
The set of decreed acting constraints $C^{\meta}(x^{\meta})$ is the set that contains all the acting constraints that are decreed by the meta component $x^{\meta}$.
\label{def:set_constraints}
\end{mydef}

Similarly to a parametrized set $\mathcal{X}^t(x^{\meta})$ defined in Definition~\ref{def:parametrized_set}, the dependency of $C^{\meta}(x^{\meta})$ with $x^{\meta}$ is defined through a parametrization with respect to to the meta component $x^{\meta}$. 

Moreover, the set of decreed acting constraints $C^{\meta}(x^{\meta})$ is a subset of the set of decreed constraints $C^{\meta}$. %, i.e. $C^{\meta}(x^{\meta}) \subseteq C^{\meta}$. 
A constraint $c \in C^{\meta}$ is either acting or nonacting, whereas $\hat{c} \in C^{\meta}(x^{\meta})$ is an acting constraint, decreed by the meta component $x^{\meta}$. In the MLP example discussed in Section~\ref{sec:example}, the set of decreed constraints is 
\begin{equation}
    C^{\meta} \ = \ \{ u_{i}- u_{i-1}\leq 0  \ : \ \forall i \in \{2,3,\ldots, l^{\max} \} \}
\end{equation}
\noindent and the set of decreed acting constraints is
\begin{equation}
   C^{\meta}(x^{\meta})=C^{\meta}(l,o)=
   \begin{cases}
   \emptyset, & \text{ if } l \in \{0,1\} \\
   \{ u_{i} - u_{i-1} \leq 0  \ : \ \forall i \in \{2,3,\ldots, l \} \} \subseteq C^{\meta}, & \text{ otherwise }
   \end{cases} 
\end{equation}
\noindent where $l \leq l^{\max}$. 

Moreover, some constraints are not decreed by the meta component $x^{\meta}$. These constraints are called the global constraints. In the MLP example, the global constraint is $c(x)= \sum_{i=1}^{l}u_i - \hat{u} \leq 0$, which is always acting no matter the meta component $x^{\meta}$.

To define the feasible set $\Omega$, the global constraints and decreed constraints are distinguished. 

\begin{mydef}[Feasible set]
\noindent The feasible set $\Omega \subseteq \mathcal{X}$ is the domain $\mathcal{X}$ defined by constraints:
\begin{align}
\begin{split}
    \Omega =  \big\{ \ (x^{\meta},x^{\cat},x^{\standard}) \in \mathcal{X} \ \ : \ \
    & c_i(x) \leq 0, \ \forall i \in \{1,2,\ldots, p \}, \\
    &c^{\meta}(x) \leq 0, \  \forall c^{\meta} \in C^{\meta}(x^{\meta})  \ \big\}
\label{def:feasible_set}
\end{split}
\end{align}
\noindent where $c_i$ are the global constraints with $p\in \mathbb{N}$ and $C^{\meta}(x^{\meta})$ is the set of decreed acting constraints, which parametrized with respect to meta component $x^{\meta}$. The number of acting constraints that are decreed by the meta component $x^{\meta}$ is simply $|C^{\meta}(x^{\meta})|$.
\end{mydef}

%--------------------------------------------------%
\subsection{Mathematical modeling of the MLP example}
\label{sec:worked_example_notation}
%--------------------------------------------------%

Each hyperparameter is identified with its variable type and role in Table~\ref{tab:hyperparameter_identified}. 

\begin{table}[htb!]
\small
\centering
\begin{tabular}{@{}llclcc@{}}
\toprule
\multicolumn{2}{l}{Hyperparameter}               & Variable & Scope     & Type   & Role  \\ \midrule
\multicolumn{2}{l}{Learning rate}                & $r$      & $]0, 1[$       & continuous & global  \\
\multicolumn{2}{l}{Activation function}          & $a$      & $ \{$ReLU, Sigmoid$\} $         & categorical & global \\
\multicolumn{2}{l}{\# of hidden layers}          & $l$      & $\{ 0,1, \ldots, l^{\max} \} $    & meta       & meta       \\
\multicolumn{2}{l}{\hspace{0.2cm} \# of units hidden layer $i$} & $u_{i}$ & $ \{u_{i}^{\min} ,u_{i}^{\min}+1, \ldots, u_{i}^{\max} \} $      & integer & decreed \\
\multicolumn{2}{l}{Optimizer}                    & $o$        & $ \{$Adam, ASGD$\} $           & meta       & meta  \\
\multicolumn{2}{l}{\hspace{0.2cm} if $o=$ ASGD}  &            &    &              \\
       & \hspace{0.4cm} decay                    & $\lambda$  & $]0, 1[$     & continuous & decreed    \\
       & \hspace{0.4cm} power update             & $\alpha$   & $]0, 1[$      & continuous & decreed    \\
       & \hspace{0.4cm} averaging start          & $t_0$      & $] 1\text{E}3, 1\text{E}8[$     & continuous & decreed    \\
\multicolumn{2}{l}{\hspace{0.2cm} if $o=$ Adam} &             &     &            &   \\
       & \hspace{0.4cm} running average 1        & $\beta_1$  & $]0, 1[$     & continuous & decreed   \\
       & \hspace{0.4cm} running average 2        & $\beta_2$  & $]0, 1[$     & continuous & decreed    \\
       & \hspace{0.4cm} numerical stability      & $ \epsilon$ & $]0, 1[$  & continuous & decreed     \\
\bottomrule
\end{tabular}
\caption{Hyperparameters with their variable type and role.}
\label{tab:hyperparameter_identified}
\end{table}

The following observations can be made. 
First, the number of units $u_i$ in the hidden layers are typed as integer variables. Although they affect the network architecture, they are not meta variables because they do not decree other variables. More precisely, they do not affect the dimension of the integer component, since they do not decree any other hyperparameters. Second, the number of hidden layers $l$ is a meta variable, since it decrees the units $u_{i}$ and thus it affects the dimension of a component $x^{\integer} \in \mathcal{X}^{\integer}(x^{\meta})$. Third, the activation function $a \in \{\text{ReLu}, \text{Sigmoid} \}$ is an unordered variable, since it is a qualitative discrete variable that belongs to a set with no appropriate metric and no order. Fourthly, the optimizer is a meta variable. Indeed, the choice of the optimizer $o$ decrees some continuous hyperparameters of the problem. 

%--------------------------------------------------%
\subsubsection{Components and sets}
%--------------------------------------------------%

The meta set $\mathcal{X}^{\meta}$ is the Cartesian product between the scopes of the two meta variables, the number of hidden layers $l$ and the optimizer $o$, thus the meta component $x^{\meta}$ and the meta set $\mathcal{X}^{\meta}$ are:

\begin{equation}
x^{\meta}  =  (l,o)
\in \mathcal{X}^{\meta}  =  \{0, 1, \ldots, l^{\max} \} \times \{\text{Adam, ASGD} \}.
\label{eq:set_meta_ex}
\end{equation}

Then, the only categorical variable is the activation function $a$, which is a global variable.  Thus, $\mathcal{X}^{\cat}(x^{\meta})=\mathcal{X}^{\cat}$ in the example, since no parametrization of the categorical set is necessary. Following this, the categorical component $x^{\cat}$ and the categorical set $\mathcal{X}^{\cat}$ are:
\begin{equation}
    x^{\cat} = a 
    \in \mathcal{X}^{\cat} =   \left \{\text{ReLU, Sigmoid} \right\}.
    \label{eq:ex_cat}
\end{equation}

Moreover, the integer component is directly the vector of units in the hidden layers, such that $x^{\integer}=u(l)=(u_1,u_2,\ldots,u_l)$. All the integer variables are decreed by the meta component $x^{\meta}$ and more specifically the number of hidden layers $l$. The integer component $x^{\integer}$ and the parametrized integer set $\mathcal{X}^{\integer}$ are:

\begin{equation}
    x^{\integer} =
    \begin{cases}
    \emptyset~\text{(nonacting)},  &\text{ if } l=0 \\
    (u_1,u_2, \ldots, u_l) \in \mathcal{X}^{\integer}(x^{\meta}) =  \mathcal{X}^{\integer}(l) = \prod_{i=1}^{l} \{u_{i}^{\min},u_{i}^{\min}+1, \ldots, u_{i}^{\max} \} \subseteq \mathbb{N}^{l}, & \text{ if } l\geq 1
    \end{cases}
\end{equation}
\noindent where $u_i^{\min}$ and $u_i^{\max}$ are respectively the minimum and the maximum of units allowed for each hidden layer $i \in \{1,2,\ldots, l \}$, $l$ is the number of hidden layers and $u(0)$ is an empty vector.

Finally, all continuous variables are decreed by the optimizer $o$, except for the learning rate $r$. Thus, the continuous component $x^{\continuous}$ is decreed by the meta component $x^{\meta}$, implying that the continuous set requires a parametrization. The continuous component $x^{\continuous}$ and the parametrized continuous set $\mathcal{X}^{\continuous}(x^{\meta})$ are:

\begin{equation}
x^{\continuous} \in \mathcal{X}^{\continuous}(x^{\meta}) =   
\begin{cases}
\mathcal{X}^{\continuous}(\text{Adam}) = ]0,1[^4 \times    \subseteq  \mathbb{R}^4 , & \text{ if $o=$ Adam } \\ 
\mathcal{X}^{\continuous}(\text{ASGD}) = ]0,1[^3 \times ] 1\text{E}3, 1\text{E}8[ \subseteq \mathbb{R}^4,  & \text{ if $o=$ ASGD }.
\end{cases}
\end{equation}

For the sake of simplicity, the scope of the units in the hidden layers $u_i$ and the number of hidden layer $l$ are set as: $u_i^{\min}=100$ and
$u_{i}^{\max}=300, \ \forall i$ and $l \in \{2, 3 \} $ in Table~\ref{tab:hyperparameter}. With $l \in \{2,3 \}$, the meta set~\eqref{eq:set_meta_ex} can be explicit as: 

\begin{equation}
\mathcal{X}^{\meta} = \{ (\text{Adam},2), (\text{Adam},3) , (\text{ASGD},2), (\text{ASGD},3)   \}.
\end{equation}

Moreover, the parametrized integer set $\mathcal{X}^{\integer}(x^{\meta})$ can also be explicit: 

\begin{equation}
\mathcal{X}^{\integer}(x^{\meta})  = \mathcal{X}^{\integer}(l) = \{100,101, \ldots, 300 \}^{l} =
\begin{cases}
\{ 100,101, \ldots, 300 \}^2, & \text{ if } l =  2 \\ 
\{ 100,101, \ldots, 300 \}^3, & \text{ if } l  = 3.
\end{cases}
\end{equation}

The parametrized categorical and continuous sets remain unchanged. 

%--------------------------------------------------%
\subsubsection{Constraints}
%--------------------------------------------------%

In the example there is a global constraint and decreed constraints. The global constraint can be easily expressed as
$c(x) \leq 0$ where
\begin{align}
\begin{split}
    c(x) = c(l) =  \sum_{i=1}^{l} u_i - \hat{u} 
    &=
    \begin{cases}
         u_1+u_2 - \hat{u}, &\text{ if } l=2 \\
         u_1+u_2+u_3 - \hat{u}, &\text{ if } l=3
    \end{cases}
\end{split}
\label{eq:example_classical_constraint}
\end{align}
The set of decreed constraint is
\begin{equation}
   C^{\meta} =  \big\{ u_{i}-u_{i-1} \leq 0  \ : \ i \in \{2,3,\ldots, l^{\max} \} \ \big \}
\end{equation}
and the set of acting decreed constraints is
\begin{equation}
   C^{\meta}(x^{\meta}) = C^{\meta}(l) =
   \begin{cases}
   \emptyset~\text{(nonacting)}, &\text{if } l<2, \\
   \big\{ u_{i}- u_{i-1} \leq 0 : \ i \in \{2,3,\ldots, l \}, \ \big \}, & \text{if } l\geq 2,
   \end{cases}
\end{equation}
\noindent which can be further detailed since $l \in \{2,3 \}$
\begin{equation}
    C^{\meta}(2) =  \{u_2 - u_1 \leq 0  \}, \quad
    C^{\meta}(3) =  \{u_3 - u_2 \leq 0, \ u_2 - u_1 \leq 0 \}.
\label{eq:example_decreed_constraints}
\end{equation}

In this particular example, the number of global constraint is $p=1$ and the number of acting constraints that decreed by the meta component is $|C(x^{\meta})|=l-1$.

%--------------------------------------------------%
\subsubsection{Visualization of the domain and the feasible set}
%--------------------------------------------------%

\begin{figure}[htb!]
\centering
  \includegraphics[width=1.0 \linewidth]{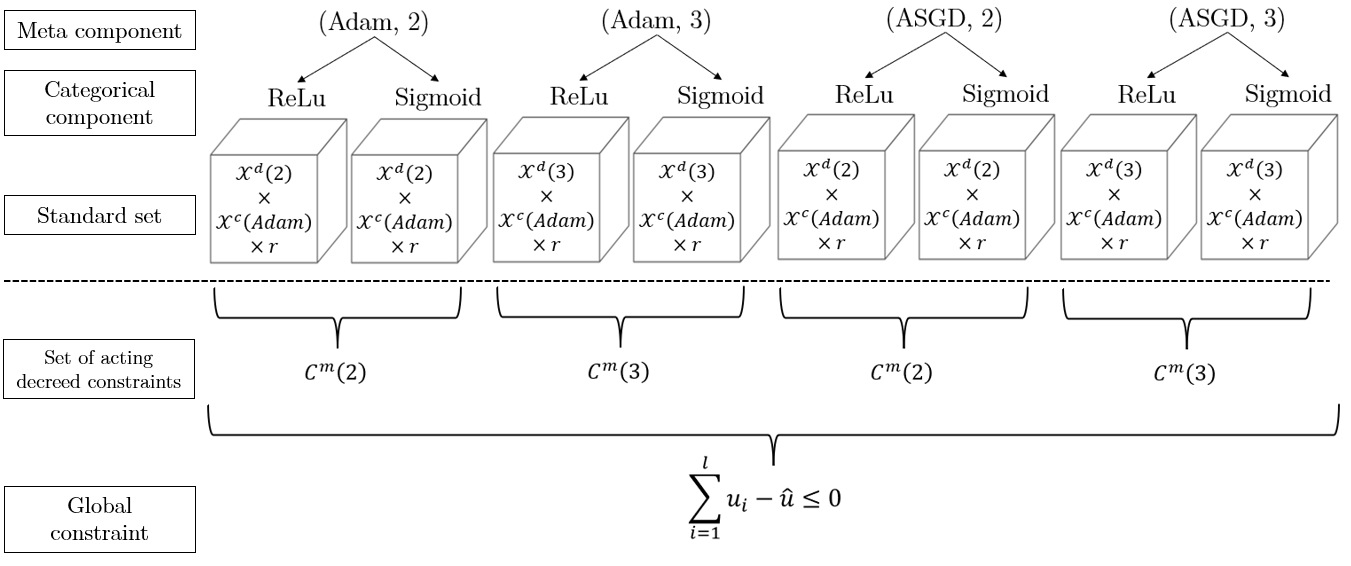}
  \caption{Diagram of the domain $\mathcal{X}$ and the constraints for the MLP example.}
  \label{fig:set_example}
\end{figure}

The alternative formulation of the domain $\mathcal{X}$ in Equation~\eqref{eq:domain2} and the feasible set $\Omega$ in Definition~\ref{def:feasible_set} of the MLP example can be visualized in Figure~\ref{fig:set_example}.

The upper part of Figure~\ref{fig:set_example}  (above the dotted line) represents the alternative formulation of the domain $\mathcal{X}$ in Equation~\eqref{eq:domain2}. The parametrized standard sets $\mathcal{X}^{\standard}(x^{\meta})$ are illustrated as small boxes and the unions from left to right in Equation~\eqref{eq:domain2} are viewed from top to bottom in Figure~\ref{fig:set_example}. Moreover, the parametrized standard sets are expressed explicitly, such that $X^{\standard}(x^{\meta})= X^{\standard}(l,o) = \mathcal{X}^{\integer}(l) \times \mathcal{X}^{\continuous}(o)$. The lower part of Figure~\ref{fig:set_example} schematizes the constraints. The acting constraints decreed by a meta component $x^{\meta}$, are contained in the set of acting decreed constraints $C^{\meta}(x^{\meta})$. The global constraint is always acting and unaffected by the meta component $x^{\meta}$, hence it is not assign to a specific meta component $x^{\meta}$ comparatively to decreed constraints: this representation shows the global aspect of global constraints. Altogether, the upper and lower parts Figure~\ref{fig:set_example} synthesize the feasible set $\Omega$ of the MLP example.

In the literature review, it has been discussed that some optimization approaches tackle categorical variables by solving many subproblems in which a categorical component $x^{\cat}$ is fixed. Indeed, in~\cite{AbAuDe2007a, AuDe2006} the MADS algorithm was applied to a continuous space where a  discrete component, which contained meta, categorical and integer variables, was fixed. This idea can be generalized to the proposed notation system. For example, assume that $x^{\meta}=(\text{Adam},2)$ and $x^{\cat}=\text{ReLU}$ are selected and fixed. Then, the objective function $f$ could then be optimized on the parametrized standard $\mathcal{X}^{\standard}(\text{Adam},2)$ with both the meta and categorical components fixed. Subproblems are further discussed in the next Section~\ref{sec:algorithmic_framework} and more particularly in Section~\ref{sec:algo_subproblems}. 

%--------------------------------------------------%
\section{Solution strategies}
\label{sec:algorithmic_framework}
%--------------------------------------------------%

Most blackbox approaches in mixed-variable optimization are built upon two strategies. 
One solution strategy consists of solving many subproblems in which some selected components are fixed.
Another strategy consists of formulating a less costly problem that selects a candidate point to be evaluated by the more costly objective function $f$. 
Some methods rely on both strategies.

For example, direct search methods~\cite{AbAuDe2007a, AuDe01a,AuDe2006, KoAuDe01a} divide the main problem into many subproblems, in which the objective function $f$ is optimized on a continuous space for a fixed discrete component $x^d$.
Bayesian optimization (BO) formulates an auxiliary problem, with a fixed acquisition function and a probabilistic surrogate, and then selects a candidate point that is subsequently evaluated by the objective function $f$. 
The methodology proposed in~\cite{PeBrBaTaGu2021} formulates many auxiliary subproblems, where each problem has a fixed dimensional component~\cite{LuPi04a, LuPiSc05a} and each subproblem has its own surrogate. 

The two strategies are respectively defined as the subproblems strategy and the auxiliary problem strategy. These strategies are the basis of the general algorithmic framework, since most algorithms that tackle mixed-variable blackbox optimization conceptually rely on solving many subproblems or on an auxiliary problem. 

The purpose of this section is to illustrate that the framework notation may be easily adapted to the main blackbox approaches in mixed-variable optimization. More precisely, direct search and heuristic approaches are discussed through the subproblems strategy in Section~\ref{sec:algo_subproblems} and the BO approach is discussed through the auxiliary problem strategy in Section~\ref{sec:algo_aux}.

%The purpose of this section is to formulate a general algorithmic framework that is compatible in the main blackbox approaches. To achieve that, some additional algorithmic terminology is discussed in Section~\ref{sec:algo_terminology} in the first place. Finally, the subproblems strategy and auxiliary problem strategy are formalized respectively in Section~\ref{sec:algo_subproblems} and in Section~\ref{sec:algo_aux}.

%\subsection{Evaluation, DoE and Cache}
%\label{sec:algo_terminology}

%An evaluation of the objective function $(x,f(x))\in \mathcal{X} \times \mathbb{R}$ consists of a couple of a point $x \in \mathcal{X}$ and its corresponding image $f(x)\in \mathbb{R}$. Then, the Cache $V_k$ is the set that contains the previous evaluations that have been completed at the beginning of an iteration $k>0$. Furthermore, a design of experiment (DoE) is often done prior any optimization process. A $\text{DoE}~=~ \{(x_{(i)},f(x_{(i)}) \}_{i=1}^D$ often serves as an initial global exploration strategy in which $D$ scattered points on the domain $\mathcal{X}$ are evaluated. Moreover, surrogates and models are often initially constructed from evaluations provided by a DoE. Finally, the set of all the completed evaluations at iteration $k$ is the union of the DoE and the Cache $V_k$, such that $\mathbb{X}_k = \text{DoE} \cup V_k$ is called the set of evaluations. Note that, the DoE is not always executed, which partly why the Cache $V_k$ and the DoE are defined separately. 

%--------------------------------------------------%
\subsection{Subproblems}
\label{sec:algo_subproblems}
%--------------------------------------------------%

The motivation of dividing a main problem into many subproblems arises from two rationales: 1)~there are  methods that treat standard problems, or even categorical-standard problems (mostly with an auxiliary problem strategy); 2)~ there are few efficient methods that address mixed variable optimization problems with both meta (or dimensional) and categorical variables. 

%A few things to note beforehand discussing the subproblems strategy. First, 
In the context of this work, subproblems are 
obtained by fixing values of meta and categorical components.
%related to meta and categorical variables, which are the variables that are bothersome. More precisely, the meta and categorical components are the components that are interesting to fix. 
In~\cite{AbAuDe2007a, AuDe01a,AuDe2006, KoAuDe01a}, the component that is fixed is the discrete component, which contains categorical variables. Secondly, note that there's no particular interest fixing the integer or continuous components, since they can be properly optimized in practice.

To further formalize the subproblems, the objective subfunction must be first defined.    

\begin{mydef}[Objective subfunction]
An objective subfunction $g$ is the objective function $f$ with a single or many fixed components. The objective subfunction is said to be parametrized with respect to to the fixed component(s).
\label{def:objective_subfunction}
\end{mydef}

From Definition~\ref{def:objective_subfunction}, it should be noted that there is a direct correspondence between the fixed component(s) and its subproblem. In other words, a specific subproblem may be referred by its fixed component(s). Again, in Definition~\ref{def:objective_subfunction}, the components that are interesting to fix are the meta component $x^{\meta}$ and the categorical component $x^{\cat}$. In this work, only the standard subproblems, in which both the meta and categorical components, is detailed. Remember that term standard encapsulates integer and continuous. 

%--------------------------------------------------%
\subsubsection{Standard subproblems}
\label{sec:standard_subproblem}
%--------------------------------------------------%

In the standard subproblems strategy, the meta component $x^{\meta}$ and the categorical component $x^{\cat}$ are fixed, in order to generate standard subproblems (one per couple $(x^{\meta}, x^{\cat})$). Fixing a meta component $x^{\meta} \in \mathcal{X}^{\meta}$  simplifies the optimization problem, since the acting variables, the acting constraints and the dimension in the subproblems are determined. In addition, fixing the categorical components also further simplifies the optimization problem. Indeed, with both the meta and categorical components fixed, the subproblems are a standard blackbox optimization problem, where the acting variables are either integer or continuous variables. In practice, there are efficient methods to tackle these standard subproblems.   

For the standard subproblems strategy, the objective subfunction  $g: \mathcal{X}^{\standard}(x^{\meta}) \to \mathbb{R}$, parametrized with respect to the meta component $x^{\meta}\in \mathcal{X}^{\meta}$ and the categorical component $x^{\cat} \in \mathcal{X}^{\cat}(x^{\meta})$, is defined as:
\begin{align}
    g(x^{\standard};x^{\cat},x^{\meta}) = f(x^{\meta},x^{\cat},x^{\standard}), \ \text{ where } x^{\meta} \in \mathcal{X}^{\meta} \text{ and } x^{\cat} \in \mathcal{X}^{\cat}(x^{\meta}) \text{ are fixed.}
\end{align}
\noindent Thus, for a fixed meta component $x^{\meta} \in \mathcal{X}^{\meta}$ and a fixed categorical component $x^{\cat} \in \mathcal{X}^{\cat}(x^{\meta})$, a standard subproblem may be formulated as

\begin{eqnarray}
\begin{array}{rcll}
    (P^{\standard})  \hspace{1cm} 
        & \displaystyle \min_{x^{\standard} \in        
           \mathcal{X}^{\standard}(x^{\meta})} 
        & g(x^{\standard};x^{\meta},x^{\cat}) \\
    &\text{s.t.} &c^{\meta}(x) \leq 0, 
    &\forall c^{\meta} \in C^{\meta}(x^{\meta}), \\
    &&c_i(x) \leq 0, 
    &\forall i \in \{1,2,\ldots, p\}.
\end{array}
\label{eq:standard_subproblem}
\end{eqnarray}
where $P^{\standard}$ stands for standard subproblem. Moreover, note that the constraints of the problem are treated directly within the subproblems of the form $(P^{\standard})$. 

%--------------------------------------------------%
\subsubsection{Exploration of subproblems}
%--------------------------------------------------%

There is a direct correspondence between the fixed component(s) and their subproblem, hence the exploration of subproblems may be done accordingly to the fixed components. Solving subproblems may be done directly with simple heuristics, such as random searches on the meta and categorical components. However, extra work is required in a direct search framework. Qualitative variables, such as the categorical variables, do not posses intuitive neighborhoods nor directions of exploration. Hence, the meta set $\mathcal{X}^{\meta}$, which may contain meta components with meta categorical variables, and the parametrized categorical set $\mathcal{X}^{\cat}$ are both endowed with a user-defined neighborhood mapping. To formalize the exploration of subproblems, the following definition based on~\cite{AbAuDe2007a, AuDe01a,AuDe2006, KoAuDe01a}, is proposed.

%Moreover, since the categorical component $x^{\cat} \in \mathcal{X}^{\cat}(x^{\meta})$ may contain decreed variables, a subtlety emerges regarding the artificial neighborhood mapping in the categorical case $t=\cat$. Hence, for a given component $x^{\cat} \in \mathcal{X}^{\cat}(x^{\meta})$, a rule $r \in \mathcal{R}^{\cat}$, the set of rules $\mathcal{R}^{\cat}$ and the artificial neighborhood mapping $\mathcal{N}^{\cat}$ are all parametrized with respect to the meta component $x^{\meta} \in \mathcal{X}^{\meta}$. More precisely, for a given categorical component $x^{\cat} \in \mathcal{X}^{\cat}(x^{\meta})$, the set of rules is $\mathcal{R}^{\cat}(x^{\cat}; x^{\meta})$, then a neighbor is $y^{\cat}=r(x^{\cat};x^{\meta})$ and the artificial neighborhood of $x^{\cat}$ is $\mathcal{N}^{\cat}(x^{\cat};x^{\meta})\subseteq \mathcal{X}^{\cat}(x^{\meta})$. 

%\todo[inline]{EHH : le user-defined neighborhood mapping a été modifié en bleu ci-dessous. J'ai ajouté une remarque pour la propriété de décret pour le cas $\mathcal{N}^{\cat}$}

\begin{mydef}[User-defined neighborhood mapping]
For any $t \in \{\meta, \cat \}$, a user-defined neighborhood mapping $\mathcal{N}^t$ assigns a user-defined neighborhood $\mathcal{N}^t(x)\subseteq \mathcal{X}^t$ to a point $x \in \mathcal{X}$, such that each neighbor $y^t \in \mathcal{N}^t(x)$ is component of type $t$ that is determined by a given rule $r^t:\mathcal{X} \to \mathcal{X}^t$:

\begin{align}
\begin{split}
    \mathcal{N}^t \ : \ \mathcal{X} &\to \mathcal{P}(\mathcal{X}^t) \\
    x &\mapsto  \ \big\{ \ y^t \in \mathcal{X}^t \ : \ y^t = r^t(x),  \  r^t \in \mathcal{R}^t(x)  \big\} \subseteq \mathcal{X}^t    
\end{split}
\label{eq:neighborhood}
\end{align}
%or equivalently,
%\begin{equation}
%     \mathcal{N}^t(x^t) = \big\{ \ y^t \in \mathcal{X}^t \ : \ y^t = r(x^t),  \  r \in \mathcal{R}^t(x^t) \big\}  \subset \mathcal{X}^t 
%\label{eq:neighborhood}
%\end{equation}

\noindent where $r^t\in \mathcal{R}^t(x)$ is a rule that assigns a neighbor $y^t = r^t(x)\in \mathcal{X}^t$ to a point $x \in \mathcal{X}$, $\mathcal{R}^t(x)$ is a set of rules defined for the given point $x\in \mathcal{X}$ and $\mathcal{P}(\mathcal{X}^t)$ is the powerset of  $\mathcal{X}^t$, which is denoted as the codomain of the mapping $\mathcal{N}^t$ to indicate $\mathcal{N}^t(x)$ can either be:
\begin{enumerate}
    \item $\mathcal{N}^t(x) = \emptyset$, such that $x$ has no neighbor of type $t$;
    \item $\mathcal{N}^t(x) = \{ y^t \}$, such that $x$ has a single neighbor of type $t$;
    \item $\mathcal{N}^t(x) \subseteq \mathcal{X}^t$, such that $x^t$ has multiple neighbors of type $t$.
\end{enumerate}
\label{def:set_neighbor_struct}
\end{mydef}

The set of rules $\mathcal{R}^t(x)$ embeds the generality of the user-defined neighborhood $\mathcal{N}^t(x)$. 
Indeed, a rule $r^t \in \mathcal{R}^t(x)$ must only respect the following mapping $r: \mathcal{X} \to \mathcal{X}^t$, which indicates that a component $y^{t}=r^t(x) \in \mathcal{X}^t$, called a neighbor, is assigned to a point $x \in \mathcal{X}$. 
In practice, it is from the these rules that user-defined neighborhoods are generated and implemented. Moreover, two issues are specific to the categorical case $t=\cat$: 1) the set $\mathcal{X}^\cat$ is the parametrized categorical set:  $\mathcal{X}^t=\mathcal{X}^{\cat}(x^{\meta})$; 2) the user-defined neighborhood mapping $\mathcal{N}^{\cat}$ takes a point $x \in \mathcal{X}$ as an argument, which allows to take into account the decree property of meta variables for the user-defined neighborhood mapping $\mathcal{N}^{\cat}$ and its constituent parts, such as the rules $r^{\cat}$.

%Moreover, since the categorical component $x^{\cat} \in \mathcal{X}^{\cat}(x^{\meta})$ may contain decreed variables, a subtlety emerges regarding the user-defined neighborhood mapping in the categorical case $t=\cat$. Hence, for a given component $x^{\cat} \in \mathcal{X}^{\cat}(x^{\meta})$, a rule $r \in \mathcal{R}^{\cat}$, the set of rules $\mathcal{R}^{\cat}$ and the user-defined neighborhood mapping $\mathcal{N}^{\cat}$ are all parametrized with respect to the meta component $x^{\meta} \in \mathcal{X}^{\meta}$. More precisely, for a given categorical component $x^{\cat} \in \mathcal{X}^{\cat}(x^{\meta})$, the set of rules is $\mathcal{R}^{\cat}(x^{\cat}; x^{\meta})$, then a neighbor is $y^{\cat}=r(x^{\cat};x^{\meta})$ and the user-defined neighborhood of $x^{\cat}$ is $\mathcal{N}^{\cat}(x^{\cat};x^{\meta})\subseteq \mathcal{X}^{\cat}(x^{\meta})$. 

In the MLP example and using Equation~\eqref{eq:set_meta_ex}, the meta rules of the form $r^{\meta}:\mathcal{X}\to \mathcal{X}^{\meta}$, for a given point $y=(y^{\meta},x^{\cat},x^{\standard})\in \mathcal{X}$ with $y^{\meta}=(l,o)$, could be
$$ \begin{array}{lll}
& r_1^{\meta}(y) = (l+1, o), & r_2^{\meta}(y) = (l-1, o), \\ r_3^{\meta}(y) = (l, \bar{o}), &
r_4^{\meta}(y) = (l+1, \bar{o}), & r_5^{\meta}(y) = (l-1, \bar{o}), \end{array}$$
where $\bar{o}$ represents the other optimizer available.
The set of rules would be:
\begin{equation}
    \mathcal{R}^{\meta}(y)= 
    \begin{cases}
        \big\{r_1^{\meta}, r_3^{\meta}, r_4^{\meta} \big\},& \text{ if } l=0 \\
        \big\{r_2^{\meta}, r_3^{\meta}, r_5^{\meta} \big\}, & \text{  if } l=l^{\max} \\
        \big\{r_1^{\meta}, r_2^{\meta}, r_3^{\meta}, r_4^{\meta}, r_5^{\meta} \big\}, & \text{ otherwise,}
    \end{cases}
\end{equation}
with corresponding user-defined neighborhood
\begin{equation}
\mathcal{N}^{\meta}(y) =
    \begin{cases}
    \big\{(l+1, o), (l, \bar{o}), (l+1, \bar{o}) \big\},  &\text{ if } l=0 \\
    \big\{(l-1, o), (l, \bar{o}), (l-1, \bar{o})   \big\}, &\text{  if } l=l^{\max} \\
    \big\{(l+1, o), (l-1, o), (l, \bar{o}), (l+1, \bar{o}), (l-1, \bar{o})   \big\},  &\text{ otherwise.}
    \end{cases}
\end{equation}

The evaluations of the blackbox objective function $f$ are generally costly, which implies that the user-defined neighborhood mappings have to set a trade-off between being exploratory and practically computational. Again, in practice, the user-defined neighborhood mappings $\mathcal{N}^{\meta}$ and $\mathcal{N}^{\cat}$ are based on rules provided by a user. Thus, the compromise is set with the discretion of the user. To lower the number of evaluations, some polling strategies may be used in practice. Indeed, instead of exploring all the neighbors at a given iteration, an opportunistic strategy would stop the iteration if a neighbor that offers a better solution is determined and resume from that neighbor.

%----------------------------------------%
\subsubsection{Direct search framework}
%----------------------------------------%

Direct search methods with strict decrease are iterative algorithms that start with an initial point $x_{(0)}$ and seek a candidate point $t$ whose objective function value $f(t)$ is strictly less than $f(x_{(k)})$, where $x_{(k)}$ is the current incumbent solution at iteration $k$. More precisely, at every iteration $k$, a set of trial points $T$ is generated. 
Opportunistically, if a trial point $t\in T$ improves the objective function value, then it becomes the next incumbent solution $x_{(k+1)}=t$ and the iteration $k$ terminates. 
Otherwise, the current incumbent solution remains unchanged, such that $x_{(k+1)}=x_{(k)}$~\cite{AuHa2017, AuLeDTr2018}. In practice, stopping the  iteration opportunistically reduces the number of evaluations required~\cite{AuHa2017}. 

Moreover, direct search methods tackle blackbox optimization problems with two main mechanism: a global search strategy (diversification) and a poll that locally searches better solutions (intensification). 

By its own, a poll is prone to miss out good point solutions. Indeed, the poll may get caught in a region with local minima or may neglect the exploration of promising regions that are far from the poll. For the meta set $\mathcal{X}^{\meta}$ and parametrized categorical set $\mathcal{X}^{\cat}(x^{\meta})$ the poll may be emulated with some user-defined neighborhood mapppings $\mathcal{N}^{\meta}$ and $\mathcal{N}^{\cat}$ respectively. The quality of a poll based on a user-defined neighborhood mapping, such as the meta and categorical polling, depends on the exhaustiveness of the set of rules $\mathcal{R}^{\meta}$ and $\mathcal{R}^{\cat}$. Therefore, depending on the quality of implementation by the user and the dimensions of the problem, the poll, based on user-defined neighborhood mappings, is likely to neglect some promising components.

In that regard, a global search may help overcome this problem by evaluating scattered trial points (or components) with a flexible strategy that serves as a diversification mechanism. The global search is generally being done before the poll for opportunistic reasons, given that the global search may find a better or interesting point that deserves to be further explored with the poll. The global search is an optional step that often improves the overall quality of a solution and increases the convergence speed. 
Many generic and low-cost global search strategies exits, such as the random search, Latin hypercube sampling
 or a Nelder-Mead search~\cite{AuTr2018},
 and more sophisticated and costly global search strategies can be implemented to generate promising trial points or unexplored regions, such as the Gaussian Processes (surrogate) paired with an acquisition function (auxiliary problem strategy) that quantifies the uncertainty and the potentiality of a point.

Algorithm~\ref{algo:direct_search} presents the main steps of a direct search methodology. 
The methodology consists of a standard subproblems strategy (see Section~\ref{sec:standard_subproblem}) paired with an exploration of subproblems that is done with user-defined neighborhood mappings $\mathcal{N}^{\meta}$ and $\mathcal{N}^{\cat}$ from Definition~\ref{eq:neighborhood}.

% Algo de base
\begin{algorithm}[htb!]

\SetAlgoLined

% Boucle principale

\While{ stopping criteria not reached }{
1. \textbf{Global search} \;
\hspace{0.2cm}  Select $t^{\meta} \in \mathcal{X}^{\meta}$ with a global meta exploration strategy \;

\hspace{0.2cm}  Select $t^{\cat} \in \mathcal{X}^{\cat}(x^{\meta})$ with a global categorical exploration strategy \;

\hspace{0.2cm} Let $t$ be obtained by solving the subproblem $\left(P^{\standard}\right)$ with $t^{\meta}$ and $t^{\cat}$ fixed

\uIf{ $f(t) < f \left(x_{(k)} \right)$ }{
$x_{\left(k+1\right)} \leftarrow t$ \;
}
\uElse{
\;
2. \textbf{Poll on user-defined neighborhoods} \;
\For{$t^{\meta} \in \mathcal{N}^{\meta} \left(x^{\meta}_{(k)} \right)$}{
        \For{$t^{\cat} \in \mathcal{N}^{\cat}\left(x^{\cat}_{(k)};t^{\meta}\right)$ }{
        Let $t$ be obtained by solving the subproblem $\left(P^{\standard}\right)$ with $t^{\meta}$ and $t^{\cat}$ fixed
        
        %Let $t$ be obtained by solving the subproblem $\displaystyle \min_{x^{\standard} \in \mathcal{X}^{\standard}(t^{\meta},t^{\cat})}~ g(x^{\standard};t^{\meta},t^{\cat})$\;

        \If{ $f(t) < f \left(x_{(k)} \right)$ }{
        $x_{(k+1)} \leftarrow t$ \;
        \FuncSty{break} \# Opportunistic strategy\; 
        } 
    }
   }
  }
} \;
\caption{Direct search main steps.}
\label{algo:direct_search}
\end{algorithm}

In Algorithm~\ref{algo:direct_search}, the two main steps to tackle the meta and categorical variables with a direct search approach are compactly presented. 
For the global search and poll steps, a standard subproblem $(P^{\standard})$, which respects the formulation in Problem~\eqref{eq:aux_problem}, is solved. Hence, the constraints of the problem are handled within the subproblems. Moreover, the solving of a subproblem $(P^{\standard})$ encapsulates many algorithmic details, such as a stopping criteria for a subproblem, as well as a global search and poll on the integer and continuous (standard) variables. Note that, a potential solver for the subproblems could be the MADS algorithm~\cite{AuDe2006} which enables to treat simultaneously integer and continuous variables (standard problem). For more details, see~\cite{AuLeDTr2018}. Then, additionally, constraints can be handle with the progressive barrier technique~\cite{AuDe09a}.

%--------------------------------------------------%
\subsection{Auxiliary problem}
\label{sec:algo_aux}
%--------------------------------------------------%

Auxiliary problems inexpensively allow to select candidate points to be evaluated by the true objective function $f$. 
Auxiliary problems are generally built from a surrogate model $\tilde{f}$ of the objective function $f$, an acquisition function $\alpha$, as well as surrogates of each global constraint $\tilde{c}_j,~j \in \{1,2,\ldots, p \}$ and decreed constraint $\tilde{c}^{\meta} \in C^{\meta}$.
The acquisition function $\alpha$ allows to select candidate points in promising regions (intensification) or in unexplored regions (exploration). 
The acquisition function $\alpha$ is generally applied to a surrogate model $\tilde{f}$ that quantifies the uncertainty of a point of its domain, and provides a prediction of the true objective function $f$. 
This is the case in BO where  $\tilde{f}$ is a GP probabilistic surrogate model .
Other surrogate models can be considered, such as random forests, however the most common remains the GPs. In this section, only GP surrogate models are adapted to the notation framework, since they are the basis of BO, an important blackbox approach to tackle mixed-variable problems. Before discussing BO, the encoding of variables is discussed.

%---------------------------------------------------------%
\subsubsection{Encoding of variables and auxiliary domain}
%---------------------------------------------------------%

BO methodologies (from Section~\ref{sec:literature_review}) often tackle categorical variables by encoding them as quantitative variables. 
For instance, the categorical variables may be encoded by the emerging latent variables or simply with the popular one-hot encoding binary vectors relaxed into a continuous vector~\cite{GaHe2020}.

%Moreover, some encoding strategies could be employed on the meta variables of any type $t\in\{ \}$. For example, latent variables could encode correlated meta components closely in a latent space, which could then serve as a basis for exploring different meta components.

\begin{mydef}[Encoder]
\noindent For any $t \in \{ \cat, \unordered, \ordered \}$ and  iteration $k \in \mathbb{N}$, the encoder $\phi^t_{(k)}$, parametrized with respect to the meta component $x^{\meta} \in \mathcal{X}^{\meta}$, is a mapping that that assigns an encoded component $l^t$ to a component $x^t$, such that
\begin{align}
    \begin{array}{lccl}
        \phi^t_{(k)} ~:& \mathcal{X}^t(x^{\meta}) &\to& \mathcal{L}^t(x^{\meta}) \\
        &x^t  &\mapsto&  l^t = \phi^t_{(k)}(x^t;x^{\meta}). 
    \end{array}
\end{align}
\end{mydef}
An encoder $\phi_{(k)}^t$ may be updated at every iteration $k\in\mathbb{N}$, such as the latent variables discussed in Section~\ref{sec:literature_review}. 
In order to take into account the decree properties of the meta variables, an encoder is parametrized with respect to the meta component $x^{\meta} \in \mathcal{X}^{\meta}$. 
In general, a meta variable may be meta categorical variable. However, in this work, the meta variables are not encoded for two reasons. First, the decreeing property of encoded meta variables may be ambiguous and difficult to conserve through sophisticated mappings, such as the latent variables. Secondly, there are categorical kernels that allow to avoid  encoding categorical variables, hence in a BO framework, meta categorical variables may be treated with these kernels.  

One of the main purpose of encoding categorical variables (or equivalently categorical component) is to formulate an auxiliary problem in which these encoded variables possess mathematical properties, making them easier to manipulate. However, by encoding the categorical variables, the domain of the surrogate model may differ from the domain of the objective function $\mathcal{X}$. Hence, the auxiliary domain $ \mathcal{X}_{\aux}$ is defined as follows.

\begin{mydef}[Auxiliary domain]
\noindent The auxiliary domain at an iteration $k \in \mathbb{N}$ is defined by:
\begin{align}
\begin{split}
    \mathcal{X}_{\aux} =  \big\{ \ \ (x^{\meta}, l^{\cat}, x^{\standard}) \ \ : \ \
    &x^{\meta} \in \mathcal{X}^{\meta}, \\
    &l^{\cat} \in \mathcal{L}^{\cat}(x^{\meta}), \\
    &x^{\standard} \in \mathcal{X}^{\standard} (x^{\meta})\ \ \big\}
\end{split}
\end{align}
where $l^{\cat}=\phi_{(k)}^{\cat}(x^{\cat};x^{\meta})$ and $\mathcal{L}^{\cat}(x^{\meta})$ is the encoded parametrized categorical set.
\label{def:aux_domain}
\end{mydef}

Definition~\ref{def:aux_domain} allows to set $\mathcal{L}^{\cat}(x^{\meta})=\mathcal{X}^{\cat}(x^{\meta})$, so that no encoding is done: $l^{\cat}=x^{\cat}$. 
In addition, since some categorical kernels do not require encoding, it follows that the auxiliary domain $\mathcal{X}_{\aux}$ is compatible with encoded categorical variables or with the original categorical variables.

From Definition~\ref{def:aux_domain}, the auxiliary maximization problem may be formulated as: 
\begin{eqnarray} \renewcommand{\arraystretch}{1.2}
\begin{array}{lcll}
    (P^{\aux})  \hspace{1cm}    
    &\displaystyle \max_{x \in \mathcal{X}_{\aux}} &  \alpha \left(x;\tilde{f} \right) \\
    &\text{s.t.} 
    &\tilde{c}_i(x) \leq 0, 
        & \forall i \in \{1,2,\ldots, p\} \\
    &&\tilde{c}^{\meta}(x) \leq 0, 
        &\forall \tilde{c}^{\meta} \in \tilde{C}^{\meta}(x^{\meta}), \\
    &&x^{\cat} =\phi_{(k)}^{\cat}(x^{\cat};x^{\meta})
        & \mbox{for some }  x^{\cat} \in \mathcal{X}^{\cat}(x^{\meta}),
\end{array}
\label{eq:aux_problem}
\end{eqnarray}

\noindent where $(P^{\aux})$ stands for auxiliary problem, $\alpha:\mathcal{X}_{\aux} \to \mathbb{R}$ is an acquisition function applied to a surrogate model $\tilde{f}$, $\tilde{c}_i \ \forall i \in \{1,2,\ldots, p \}$ are surrogate constraints for the global constraints, $\tilde{c}^{\meta} \in \tilde{C}^{\meta}$ is a surrogate constraint for a decreed constraint. The last constraint imposes the existence of some $x^{\cat} \in \mathcal{X}^{\cat}(x^{\meta})$ such that $x^{\cat} =\phi_{(k)}^{\cat}(x^{\cat};x^{\meta})$ is a pre-image constraint that recovers a categorical component $x^{\cat} \in \mathcal{X}^{\cat}(x^{\meta})$ from the encoded parametrized categorical set $\mathcal{L}^{\cat}(x^{\meta})$.
The pre-image constraint also ensures that the optimal auxiliary problem solution resides in the domain $\mathcal{X}$.
For more details on pre-images problem, refer to~\cite{CuLeRoPeDuGl2021}. 
%{\bl The $\argmax$ in $(P^{\aux})$ outlines that a candidate point is selected by the auxiliary problem, which is more commonly a maximisation problem.  }

%There exist many surrogate for the global or decreed constraints. A popular one is the the probability of feasibility, such that a given constraint $c_j(x)\leq 0$ its surrogate constraint is of the form $\tilde{c}_j(x)=\mathbb{P}[c_j(x) \leq 0] \in [0,1]$. The surrogate constraints maybe be grouped together, such that
%\begin{equation}
%    P \left[x \text{ is feasible} \right] \ = \ \prod_{i=1}^p P \left[c_i(x) \leq 0 \right] \prod_{c_j^{\meta} \in C^{\meta}(x^{\meta})} P \left[ c_j^{\meta}(x) \leq 0 \right],
%\end{equation}
%where $p$ is the number of global constraints and  $C^{\meta}(x^{\meta})$ is the set of decreed acting constraints. 

%---------------------------------------------------------%
\subsubsection{Bayesian optimization}
\label{sec:BO}
%---------------------------------------------------------%

In this section, the BO approach is formulated as an auxiliary problem $(P^{\aux})$, without detailing the algorithmic steps or the construction of the GP (see~\cite{RaWi06} or~\cite{SSWAF2015} for more details on this subject). 
For the purpose of this work, it is sufficient to formulate the BO approach as an auxiliary problem $(P^{\aux})$ and to develop the kernel from the notation framework, since the kernel almost entirely characterizes the probabilistic surrogate (GP).
A kernel $k:\mathcal{X}_{\aux} \times \mathcal{X}_{\aux} \to \mathbb{R}$ is a positive semi-definite covariance function. 
Conceptually, the kernel establishes the mathematical properties of the GP, such as the degree smoothness. 
%In future work, the BO approach will be more thoroughly developed within the notation framework developed in this work. 
    
In its simplest noise free form, a probabilistic BO distribution is built from a GP, which allows to compute for any given point $x \in \mathcal{X}_{\aux}$, a prediction $\hat{f}(x)$ and an uncertainty measure $\hat{\sigma}^2(x)$, such that
\begin{equation}
\begin{cases}
    \hat{f}(x)  &=  \ \kappa^{\intercal}(x)K^{-1}  f\left( \mathbb{X} \right) \\
    \hat{\sigma}(x)^2 &= \  k(x,x) -  \kappa^{\intercal}(x)K^{-1}\kappa(x)
    \end{cases}
\label{eq:BO_pred_var}
\end{equation}

\noindent where $\mathbb{X}$ is a set of sample points, $f\left(\mathbb{X}\right)$ is the vector of objective function values of the sample points, $\kappa(x)$ is a vector in which an element is the computed kernel $k(x,y)$ with $(x, y) \in \mathcal{X}_{\aux} \times \mathbb{X}$, $K$ is matrix in containing all pairs $(y,z) \in \mathbb{X} \times \mathbb{X}$, such that an element of $K$ is $k(y,z)$. In~\eqref{eq:BO_pred_var}, everything is computed from the kernel $k$. In other words, Equation~\eqref{eq:BO_pred_var} displays that the GP is entirely characterized by the kernel: it is assumed that the GP is noise free and that the mean function is zero, which is a common practice~\cite{RaWi06}. 
The surrogate probabilistic model $\tilde{f}$ satisfies
\begin{equation}
 \tilde{f}(x) \sim \mathcal{N} \left( \hat{f}(x), \hat{\sigma}(x)^2 \right).  
 \label{eq:GP}
\end{equation}
where $\mathcal{N}$ is the normal distribution. Moreover, a common acquisition function $\alpha$ applied on GP surrogates is the $EI$ from~\cite{JoScWe1998}:
\begin{equation}
    EI\left(x;\tilde{f} \right)  =  \mathbb{E} \left[ \max(f_{\star}-\tilde{f}(x), 0)  \right] = \left( f_{\star} - \hat{f}(x) \right) \Phi \left( \frac{f_{\star} - \hat{f}(x)}{\hat{\sigma}(x)} \right) + \hat{\sigma}(x) \phi \left( \frac{f_{\star} - \hat{f}(x)}{\hat{\sigma}(x)} \right)
    \label{eq:EI}
\end{equation}
where $f_{\star}=f(x_{\star})$ is current best known objective function value,  $\hat{\sigma}(x)$ is the standard deviation of the GP, $\Phi$ and $\phi$ are respectively the cumulative distribution and the density function of a standard normal distribution (centered at zero with variance of one). In Equation~\eqref{eq:EI}, the intensification and exploration trade-off of the $EI$ (acquisition function) is displayed by the two terms: the first term favors promising low surrogate values (intensification) and the second term favors highly uncertain points (exploration). In the auxiliary problem $(P^{\aux})$, the acquisition function could be $\alpha (x; \tilde{f})= EI(x;\tilde{f})$.

In a similar manner to surrogate model $\tilde{f}$ evaluated at a point $x \in \mathcal{X}_{\aux}$ in Equation~\eqref{eq:GP}, the surrogate constraints in the auxiliary problem $(P^{\aux})$, may be developed into GP probabilistic surrogates. Thus, a given surrogate constraint $\tilde{c}_i$ would have its own prediction function $\hat{c}_i$ (similarly to $\hat{f}(x)$ in Equation~\eqref{eq:BO_pred_var}), which could be directly used in the auxiliary problem $(P^{\aux})$, i.e., $\tilde{c}_i(x) = \hat{c}_i(x)$. Acquisition functions may also be applied to probabilistic surrogate constraints, which is not covered in this work.
   
At this stage, the BO framework is formulated in a general manner, which does not explicit the mixed-nature of the optimization problems at stake. To adapt the BO framework on a mixed-variable context, the kernel $k$, must be further detailed with the support of the notation framework. Many possible kernels can be built with operations of multiplication and additions that respects the RKHS formalism~\cite{PeBrBaTaGu2021, RoPaDeClPeGiWy2020}. 
An example of a specific kernel is detailed next
 to illustrate the compatibility of the framework 
 with the mixed-variable optimization BO literature.

The kernel $k$ is built piece-by-piece with the partition of a point $x=(x^{\meta},x^{\unordered}, x^{\ordered}, x^{\integer},x^{\continuous})$. 
The parametrized continuous kernel $k^{\continuous}:\mathcal{X}^{\continuous}(x^{\meta}) \times \mathcal{X}^{\continuous}(x^{\meta}) \to \mathbb{R}$ is formulated as multiplication of one-dimensional squared-exponential kernels: 
    \begin{equation}
        k^{\continuous}(x^{\continuous},y^{\continuous};x^{\meta}) = \exp \left( - \sum_{i=1}^{n^{\continuous}(x^{\meta})} \lambda^{\continuous}_i \left[ x_i^{\continuous} - y_i^{\continuous} \right]^2 \right).
        \label{eq:continuous_kernel}
    \end{equation}
where the $\lambda_i^{\continuous}$ are weight coefficients (hyperparameters of the surrogate model) that can be adjusted by various methods, such as the MLE. 
%The hyperparameters of the surrogate model are not furthered discussed in this work.

% Following, based on~\cite{GaHe2020}, an approach to build a parametrized integer $k^{\integer}$ would rely on relaxed integer variables in the auxiliary domain $\mathcal{X}_{\aux}$ and apply a parametrized continuous kernel $k^{\continuous}$ with a transformation $T$ that rounds the relaxed integer variables to the closet integer:
    
The parametrized integer kernel $k^{\integer}: \mathcal{X}^{\integer}(x^{\meta}) \times \mathcal{X}^{\integer}(x^{\meta}) \to \mathbb{R}$ is similar to $k^{\continuous}$, but applies a transformation $T$ that rounds the relaxed integer variables to the nearest integer~\cite{GaHe2020}:
\begin{equation}
    k^{\integer}(x^{\integer},y^{\integer};x^{\meta})   = \exp \left( - \sum_{i=1}^{n^{\integer}(x^{\meta})} \lambda^{\integer}_i \left[ T \left(x_i^{\integer}\right) - T \left(y_i^{\integer} \right) \right]^2 \right)
    \label{eq:discrete_kernel}
\end{equation}
% k^{\continuous} \left( T(x^{\integer}),T(y^{\integer}); x^{\meta} \right)

\noindent where $x_i^{\integer},y_i^{\integer}~\forall i \in I^{\integer}(x^{\meta})$ are relaxed integer variables and $\lambda_i^{\integer}$ are hyperparameters of the surrogate model. 
The transformation $T$ conserves the order of an integer variable and ensures that the one-dimensional kernels in~\eqref{eq:discrete_kernel} are piecewise functions~\cite{GaHe2020}.

The parametrized standard kernel $k^{\standard}:\mathcal{X}^{\standard}(x^{\meta}) \times \mathcal{X}^{\standard}(x^{\meta}) \to \mathbb{R} $ is formulated as multiplication of $k^{\integer}$ and $k^{\continuous}$:
\begin{equation}
    k^{\standard}(x^{\standard},y^{\standard};x^{\meta}) = k^{\integer}(x^{\integer},y^{\standard};x^{\meta}) \cdot k^{\continuous}(x^{\standard},y^{\standard};x^{\meta}).
    \label{eq:kernel_standard}
\end{equation}

The parametrized categorical kernel $k^{\cat}$ may be formulated with an encoding on the categorical variables~\cite{ZhTaChAp2020} ($k^{\cat}:\mathcal{L}^{\cat}(x^{\meta}) \times \mathcal{L}^{\cat}(x^{\meta})  \to \mathbb{R}$), or without any encoding ($k^{\cat}:\mathcal{X}^{\cat}(x^{\meta}) \times \mathcal{X}^{\cat}(x^{\meta}) \to \mathbb{R}$). With an encoding, the parametrized categorical kernel $k^{\cat}$ is similar to $k^{\continuous}$: 

\begin{align}
    k^{\cat}(l^{\cat},u^{\cat};x^{\meta}) = \exp \Bigg( 
    \displaystyle - \sum_{i \in \mathcal{L}_{\aux}^{\cat}(x^{\meta}) } \lambda^{\cat}_i \left[ l_i^{\cat} - y_i^{\cat} \right]^2 
    \Bigg),
    \label{eq:kernel_cat1}
\end{align}
\noindent where $\lambda^{\cat}_i$ are hyperparameters of the surrogate model and $\mathcal{L}_{\aux}^{\cat}(x^{\meta})$ is the set of indices of the encoded (acting) categorical variables. Without any encoding, the parametrized categorical kernel $k^{\cat}$ is formulated as tensor products of matrices (one matrix per categorical categorical)~\cite{RoPaDeClPeGiWy2020}:
\begin{align}
    k^{\cat}(x^{\cat},y^{\cat};x^{\meta})  = \left(  \kernel{i=1}{n^{\unordered}(x^{\meta}) }  T_i^{\unordered}  \left( x_i^{\unordered}, y_i^{\unordered} \right) \right) \otimes \left( \kernel{i=j}{n^{\ordered}(x^{\meta}) } T_j^{\ordered} \left(x_j^{\ordered}, y_j^{\ordered} \right) \right),
    \label{eq:kernel_cat2}
\end{align}
\noindent where, for $t \in \{\unordered, \ordered \}$ and a categorical variable $x_i^t \in \{1,2,\ldots, c_i \}$, $T_i^{t} \in \mathbb{R}^{c_i \times c_i}$ is a positive semi-definite matrix in which an element is the correlation between two categories of the variable $x_i^t$. Hence, for two given variables with specific categories $x_i^t = c_1$ and $y_i^t=c_2$, $T_i^{t} \left( x_i^{t}, y_i^{t} \right)$ is a correlation measure between the categories $c_1$ and $c_2$.  In~Equation~\eqref{eq:kernel_cat2}, the matrices for the nominal and ordinal variables $T_i^{\unordered}$ and $T_j^{\ordered}$  are distinguished, since there exist more sophisticated matrices for the ordinal variables~\cite{RoPaDeClPeGiWy2020}. 
    
Finally, a mixed kernel $k:\mathcal{X}_{\aux} \times \mathcal{X}_{\aux} \to \mathbb{R}$, based on~\cite{PeBrBaTaGu2021}, is formulated as: 
% EHH : ici j'utilise \Big car le \left( \right) n'ajuste pas comme je le voudrais
\begin{equation}
    k(x,y) =
    \begin{cases}
    \prod_{i=1}^{n^{\meta}} k_{i}^{\meta}(x_i^{\meta},y_i^{\meta}), &\text{if}~x^{\meta}\neq y^{\meta},  \\
    \prod_{i=1}^{n^{\meta}} \Big( k_{i}^{\meta}(x_i^{\meta},y_i^{\meta}) \cdot \big[  k^{\cat}(l^{\cat},u^{\cat};x^{\meta})  k^{\standard}(x^{\standard},y^{\standard};x^{\meta})  \big]    \Big), &\text{otherwise},
    \end{cases}
\label{eq:mixed_kernel}
\end{equation}
where $k_i^{\meta}:S_i^{\meta} \times S_i^{\meta} \to \mathbb{R}$ is a one-dimensional kernel for a meta variable $x_i^{\meta}\in S_i^{\meta}$, $k^{\cat}$ is the parametrized categorical kernel that may take the form in Equation~\eqref{eq:kernel_cat1} or Equation~\eqref{eq:kernel_cat2},
$k^{\standard}$ is the parametrized standard kernel in Equation~\eqref{eq:kernel_standard}. 
In Equation~\eqref{eq:mixed_kernel}, the meta kernel $k^{\meta}:\mathcal{X}^{\meta} \times \mathcal{X}^{\meta} \to \mathbb{R}$ is implicitly decomposed into one-dimensional kernels (one per meta variable), which is again, common practice in the literature. Moreover, in Equation~\eqref{eq:mixed_kernel}, the kernel computations for the categorical and standard variables are only being done if the two points in share the same meta component: for $t \in \{\cat, \standard \}$, the kernel computation $k(x^t,y^t;x^{\meta})$ in Equation~\eqref{eq:mixed_kernel} is only done if $x^{\meta}=y^{\meta}$, which implies that $x^t$ and $y^t$ must both reside in the same parametrized set $\mathcal{X}^t(x^{\meta})$.

%---------------------------------------------------------%
\section{Conclusion}
%---------------------------------------------------------%

This work proposes a thorough notation framework for mixed-variable optimization problems. 
The framework formally models mixed-variable problems with a careful emphasis on meta and categorical variables. 
More precisely, a point $x=(x^{\meta},x^{\cat},x^{\standard})$, the domain $\mathcal{X}$ of the objective function, and the feasible set $\Omega$, are rigorously defined. 
%To our knowledge, such rigorous definitions have not been proposed for mixed-variable with meta (or dimensional in the literature) and categorical variables. 
Definitions are developed to shed the light on the intrinsic difficulties resulting from the presence of meta variables. 
Notably, for $t\in \{\cat, \standard, \unordered, \ordered, \integer, \continuous \}$, a parametrized set $\mathcal{X}^t(x^{\meta})$ elucidates that some variables of type $t$ may be acting or nonacting depending on the meta variables. 

In addition, the parametrized categorical set $\mathcal{X}^{\cat}(x^{\meta})$ and the parametrized standard set $\mathcal{X}^{\standard}(x^{\meta})$ are building blocks of the domain $\mathcal{X}$ that has two equivalent formulations, respectively in Definition~\ref{def:domain} and in Equation~\eqref{eq:domain2}. 
Both formulations provide a different perspective on mixed-variable problems.
Furthermore, the constraints are split into global decreed constraints, which allow to formulate a clear feasible set $\Omega$ in Definition~\ref{def:feasible_set}. 

In Section~\ref{sec:algorithmic_framework}, the subproblems strategy and auxiliary problem strategy are exhaustively discussed from the notation framework. These strategies allow to formally adapt the notation framework to the direct search approach through the subproblems strategy, as well as to the Bayesian optimization approach through the auxiliary problem strategy. 
%Heuristics are briefly discussed within the subproblems strategy. 
Thereby, the notation framework is shown to be compatible with most of the approaches of the literature on mixed-variable optimization with meta (or dimensional) and categorical variables.

Computational experiments will be carried out in future studies with the mathematical framework of this work as a foundation. 
Moreover, the Bayesian optimization approach will also be extensively developed, with aim of bridging the communities in optimization and machine learning. \\

%In future work, the BO approach will be more thoroughly developed within the notation framework developed in this work. 

\noindent
{\bf Data availability statement:}
No relevant for the present theoretical study.\\

\noindent
{\bf Conflict of interest statement:}
On behalf of all authors, the corresponding author states that there is no conflict of interest.

%----------------------------------%
\bibliographystyle{plain}
\bibliography{bibliography.bib}
\pdfbookmark[1]{References}{sec-refs}
%----------------------------------%

\end{document}